\documentclass[preprint,12pt]{elsarticle}
\usepackage{amsfonts}
\usepackage{amsmath}
\usepackage{mathrsfs}
\usepackage[applemac]{inputenc}
\usepackage{amssymb}
\usepackage{times, ulem, amsfonts}
\usepackage{mathpazo, amsmath, amssymb, amsthm,color}
\makeatletter
\def\ps@pprintTitle{%
   \let\@oddhead\@empty
   \let\@evenhead\@empty
   \let\@oddfoot\@empty
   \let\@evenfoot\@oddfoot
}
\makeatother
\def\q{{\mathcal Q}}
\def\v{{\mathcal V}}
\def\x{{\mathcal X}}
\def\y{{\mathcal Y}}
\def\z{{\mathcal Z}}
\def\w{{\mathcal W}}
\def\p{{\mathcal P}}
\newcommand{\R}{{\mathbb{R}}}
\newcommand{\Z}{{\mathbb{Z}}}

\newcommand\diverg{\mathop{\mbox{\rm div}}}
\def\ddj{\dot \Delta_j}

\newtheorem{theorem}{Theorem}[section]
\newtheorem{lemma}[theorem]{Lemma}
\newtheorem{definition}[theorem]{Definition}

\newtheorem{proposition}[theorem]{Proposition}
\newtheorem{remark}[theorem]{Remark}
\numberwithin{equation}{section}

\begin{document}
\title{ 
 Global large solutions and incompressible limit for the compressible Navier-Stokes equations}

 \author{Zhi-Min Chen}
 \author{Xiaoping Zhai\corref{cor1}}
 \address{School  of Mathematics and Statistics, Shenzhen University, Shenzhen 518060, China}
\cortext[cor1]{Corresponding author. Email address: zhaixp@szu.edu.cn }

\baselineskip=24pt

\begin{abstract}
The present paper is dedicated to the  global large solutions and incompressible limit for the
compressible Navier-Stokes system in $\mathbb{R}^d$ with $d\ge 2$. Motivated  by the  $L^2$ work of  Danchin and Mucha [Adv. Math. 320, 904--925, 2017] in critical Besov spaces, we  extend the solution space into an $L^p$ framework. The result implies
the existence of global large solutions initially from   large highly oscillating   velocity fields.
\end{abstract}

\maketitle


\section{Introduction}
In this paper, we  study the global well-posedness of the compressible Navier-Stokes equations in the following form:
\begin{eqnarray}\label{m}
\left\{\begin{aligned}
&\partial_t \rho + \diverg(\rho v) = 0\, ,\\
 &     \partial_t ( \rho v ) + \diverg ( \rho v \otimes v ) -\mu\Delta v-(\mu+\lambda)\nabla\diverg v+ \nabla P( \rho) = 0,\\
    &(\rho,v)|_{t=0}=(\rho_0,v_0),
\end{aligned}\right.
\end{eqnarray}
where $\rho$ is density,  $v$ is   velocity,  $\mu$  is  shear viscosity coefficient and $\lambda$ is volume viscosity coefficient. Here $\mu$ and $\lambda$ are subject to  the standard  strong parabolicity assumption:
\begin{align*}
\mu>0\quad\hbox{and}\quad
\nu:=\lambda+2\mu>0.
\end{align*}
The  pressure  $P=P(\rho)$ is  smooth function such that  $P'>0$ and
that $P(\bar\rho)=0$  for some  positive constant reference density $\bar\rho.$

As one of the most popular fluid motion model in the field of  analysis and applications,  the compressible Navier-Stokes equations system  has attracted much attention and there is  a large  literature important to mathematical analysis and fluid mechanics.
The local well-posedness for the system \eqref{m}
was proved by Nash \cite{nash} for the smooth initial data being  away from vacuum. The existence of  global smooth solutions was  obtained by Matsumura and Nishida  \cite{mat},
when  the initial data is close to the equilibrium in $H^3(\R^3)\times H^3(\R^3)$.
 In general, whether  a smooth solution
blows up in finite time is an open problem.

The global existence of weak solutions was proved by Hoff \cite{hoff1995, hoff2005} assuming  discontinuous initial data with small energy.
The global existence of large  weak solution was established  by Lions \cite{lions+1}
under the isentropic assumption, i.e. $P=A\rho^\gamma$ for $\gamma\ge \frac 95$. This $\gamma$ restriction domain  was enlarged  by Feireisl {\it et al.} 
\cite{feireisl} to that of   $\gamma>\frac 32$.
Motivated by Hoff \cite{hoff1995,hoff2005},
 Huang {\it et al.} 
 \cite{huangxiangdi} obtained the existence of  global strong solutions  with  small energy.
However, the question of the regularity and uniqueness of
weak solution is generally  open even in the case of two dimensional space.

As given by Fujita and Kato \cite{fujita},
   the classical incompressible Navier-Stokes equations  has the scaling invariance property and gives rise to  critical spaces.
This observation was introduced to the compressible Navier-Stokes equations by Danchin \cite{danchin2000, danchin2001,  danchin2007} with respect to the scaling transformation
\begin{gather}\label{scaling}
\begin{cases}
(\rho_0, v_0)\rightarrow (\rho_0(\ell x), \ell v_0(\ell x)),\\
(\rho(t,x), v(t,x))\rightarrow (\rho(\ell^2t,\ell x), \ell v(\ell^2t,\ell x)), \quad \ell>0
\end{cases}
\end{gather}
 by neglecting the pressure term $P=P(\rho)$.
Here a function space being  critical with respect to (\ref{m})  means that the norm of the space is invariant with respect to the scaling transformation (\ref{scaling}).
For example,  the product space $\dot{B}_{p,1}^{\frac {d}{p}}({\mathbb R} ^d)\times \dot{B}_{q,1}^{-1+\frac {d}{q}}({\mathbb R} ^d) $, $1\le p, q\le \infty$, is critical  for the system (\ref{m}).

In the  critical space framework, a breakthrough was made by Danchin \cite{danchin2000}, showing
  the local well-posedness of (\ref{m})
for the initial data $(\rho_0-\overline{\rho}, v_0)$ in the critical Besov space $\dot B^{\frac  d2}_{2,1}(\R^d)\times \dot B^{-1+\frac  d2}_{2,1}(\R^d)$
and the global existence   of strong solutions initially  in the vicinity of   an equilibrium in the space  $(\dot B^{\frac  d2}_{2,1}(\R^d)\cap\dot B^{-1+\frac  d2}_{2,1}(\R^d))\times \dot B^{-1+\frac  d2}_{2,1}(\R^d)$.
Inspired by  Danchin  \cite{danchin2000}, Charve and Danchin \cite{charve} and
   Chen {\it et al.} \cite{chenqionglei} obtained the global well-posedness of (\ref{m}) in the critical $L^p$
framework. The critical Besov space, used by  Charve and Danchin \cite{charve} and
   Chen {\it et al.} \cite{chenqionglei}, seems the largest one in which the system (\ref{m}) is well-posed.
Indeed, Chen {\it et al.} \cite{chenqionglei2} proved the ill-posedness of (\ref{m})
in $\dot B^{\frac  3p}_{p,1}(\R^d)\times \dot B^{\frac 3p-1}_{p,1}(\R^d)$ for $p>6$.  An alternative proof to the  results of \cite{charve,chenqionglei}
 was further obtained by Haspot  \cite{haspot}  by using the viscous effective flux. Moreover,  Danchin and He \cite{helingbing} generalized the previous results by allowing the incompressible part of the velocity in the space $\dot B^{\frac  3p-1}_{p,1}(\R^3)$ with $p\in[2,4]$.
Interested readers may also refer to \cite{ xujiang, dejardins, FZZ1, FZZ, feireisl2,huangjingchi2, lions, villani} for   stability, decay estimate and zero Mach number limit of system \eqref{m}.

\subsection{The main result and its motivation}

Recently, Danchin and Mucha \cite{danchin2018} obtained
the global existence of  regular solutions to system \eqref{m}
 with arbitrary large initial velocity $v_0$,  almost  constant
density $\rho_0$, and   large volume  viscosity $\lambda$. This  result
 strongly  relies
on the fact that the limit velocity for $\lambda\to+\infty$ satisfies the incompressible Navier-Stokes equations:
\begin{eqnarray}\label{classns}
\left\{\begin{aligned}
&V_t + V\cdot \nabla V - \mu \Delta V+\nabla\Pi=0,\\
 &    \diverg V=0,\\
    & V|_{t=0} =\p v_0,
\end{aligned}\right.
\end{eqnarray}
with the Leray projection $\p=\mathcal{I}-\q$ with $\q = \nabla\Delta^{-1}\diverg$.

More precisely, the present study is motivated by the following result in $\R^2$:
\begin{theorem}(Danchin and Mucha  \cite{danchin2018})\label{yiyang}
 Let  $\nu\ge \mu$,  $v_0 \in \dot B^0_{2,1}(\R^2)$ and  $a_0:=\rho_0-1\in\dot B^0_{2,1}(\R^2)\cap \dot B^1_{2,1}(\R^2)$
such that
$$
Ce^{C(\widetilde M+\widetilde{M}^2)}\bigl( \|a_0\|_{\dot B^0_{2,1}}+\nu\|a_0\|_{\dot B^1_{2,1}}+\|\q v_0\|_{\dot B^0_{2,1}}+ \widetilde{M}^2
 +\mu^2 \bigr)\leq \sqrt{\mu\nu}
 $$
  for  a large constant $C$ and
 \begin{equation}
 \widetilde{M}= C\|\p v_0\|_{\dot B^0_{2,1}}\exp\Bigl(\frac  C{\mu^4}\|\p v_0\|_{L^2}^4\Bigr).\nonumber
 \end{equation}
Then there exists a unique global  regular solution $(\rho,v)$ to \eqref{m} such that
\begin{align*}\label{eq:reg2D}
  & v \in C([0,\infty);\dot B^0_{2,1}(\R^2)),\qquad v_t,\nabla^2v\in L^1(0,\infty;\dot B^0_{2,1}(\R^2)),\\
 & a:=\rho -1 \in C([0,\infty); \dot B^0_{2,1}(\R^2) \cap \dot B^1_{2,1}(\R^2)) \cap
  L^2([0,\infty);\dot B^1_{2,1}(\R^2)).
\end{align*}
In addition, there holds the following estimate
\begin{align*}
&\|\q v\|_{L^\infty(0,\infty;\dot B^0_{2,1})}
    + \|a\|_{L^\infty(0,\infty;\dot B^0_{2,1})} +
      \nu\|a\|_{L^\infty(0,\infty;\dot B^1_{2,1})}\nonumber\\
    &\quad\leq   Ce^{C(M+M^2)}\bigl(\|a_0\|_{\dot B^0_{2,1}} +\nu \|a_0\|_{\dot B^1_{2,1}}+\|\q v_0\|_{\dot B^0_{2,1}} +\widetilde{M}^2+\mu^2 \bigr).
\end{align*}
\end{theorem}

As a byproduct, Danchin and Mucha  \cite{danchin2018} obtained  the convergence
 $(\rho,v)\to(1,V)$ at the order of $\nu^{-\frac 12}.$
 This result  in the high-dimensional case $d\geq3$ was additionally  given in \cite{danchin2018} under the condition that the incompressible Navier-Stokes equations (\ref{classns}) admit a global large regular solution.
Moreover, they  \cite{danchin2018}
predicted that those results in $L^2$ Besov spaces be improved to the critical framework of  $L^p$ Besov spaces. The purpose of  the present paper is  to give a positive answer to   the  prediction.

For stating our main result,  a homogeneous tempered distribution  $z= \sum_{j \in \Z}\dot \Delta_j z \in\mathcal{S}'(\R^d)$ is truncated by lower and higher oscillation parts in the following sense:
\begin{align}
z^\ell:=\sum_{2^j\nu\leq 1}\ddj z\quad\hbox{and}\quad
z^h:=\sum_{2^   j\nu>1}\ddj z \ \  \mbox{ for } \ \ z\in\mathcal{S}'(\R^d).
\end{align}
Sometimes, for convenience, we will use the notation:
\begin{align}
\|z\|^{\ell}_{\dot B^{s}_{p,1}}:=  \|z^{\ell}\|_{\dot B^{s}_{p,1}}
\ \hbox{ and }\   \|z\|^{h}_{\dot B^{s}_{p,1}}:=  \|z^{h}\|_{\dot B^{s}_{p,1}}.
\end{align}

The main result of the present paper reads:
\begin{theorem}\label{dingli}
Let $2\leq p \leq \min\{4,\, 2d/(d-2)\}$ for $d>2$,  and  $2\leq p<4$ for $d=2$.
  Assume   $a_0^\ell \in \dot B^{-1+\frac {d}{2}}_{2,1}(\R^d)$, $a_0^h \in \dot B^{\frac {d}{p}}_{p,1}(\R^d)$, $\p v_0 \in \dot B^{-1+\frac {d}{p}}_{p,1}(\R^d)$, $\q v_0^\ell \in \dot B^{-1+\frac {d}{2}}_{2,1}(\R^d)$ and $\q v_0^h \in \dot B^{-1+\frac {d}{p}}_{p,1}(\R^d)$.
 Suppose that \eqref{classns} admits  a unique global solution $$V\in C([0,\infty);\dot B^{-1+\frac  dp}_{p,1}(\R^d))\cap{L^1(0,\infty;\dot B^{-1+\frac  dp}_{p,1}(\R^d))}.$$
 Denote
$$
M:= \|V\|_{L^\infty(\R^+;\dot B^{-1+\frac  dp}_{p,1})}+\mu\|V\|_{L^1(\R^+;\dot B^{1+\frac  dp}_{p,1})}+\|V_t\|_{L^1(\R^+;\dot B^{-1+\frac  dp}_{p,1})}.
$$
 Assume $\nu\ge \mu$ and the existence of  a (large) generic  constant $C$ such that
  \begin{align}
&\big\|a_0^\ell\big\|_{\dot{B}_{2,1}^{-1+\frac {d}{2}} }+\nu\big\|  a_0^\ell\big\|_{\dot{B}_{2,1}^{\frac {d}{2}} }+\nu\| a_0^h\|_{\dot B^{\frac  dp}_{p,1}}+\big\| \q v_0^\ell\big\|_{\dot{B}_{2,1}^{-1+\frac {d}{2}} }+\|\q v_0^h\|_{\dot B^{-1+\frac  dp}_{p,1}}+M^2+\mu^2\nonumber\\
&\quad\leq C\sqrt{\mu\nu}\exp\big(-C(M+M^2)\big).
 \end{align}

Then there exists a unique global regular solution $(\rho,v)$ to \eqref{m} such that
 \begin{align*}
&\p v \in C(\R^+;\dot B^{-1+\frac  dp}_{p,1})\cap L^1(\R^+;\dot B^{1+\frac  dp}_{p,1}),\\
& a^\ell \in C(\R^+;\dot B^{-1+\frac {d}{2}}_{2,1})\cap L^1(\R^+;\dot B^{1+\frac  d2}_{2,1}),\quad
 a^h \in C(\R^+; \dot B^{\frac {d}{p}}_{p,1})\cap
  L^1(\R^+;\dot B^{\frac  dp}_{p,1}),\\
   &\q v^\ell \in C(\R^+;\dot B^{-1+\frac  d2}_{2,1})\cap L^1(\R^+;\dot B^{1+\frac  d2}_{2,1}),\quad
    \q v^h \in C(\R^+;\dot B^{-1+\frac  dp}_{p,1})\cap L^1(\R^+;\dot B^{1+\frac  dp}_{p,1}).
 \end{align*}
In addition, the following estimate holds true:
\begin{align*}
&\big\|a^\ell\big\|_{\widetilde{L}^\infty(0,\infty;\dot{B}_{2,1}^{-1+\frac {d}{2}} )}
+\nu\big\|  a^\ell\big\|_{\widetilde{L}^\infty(0,\infty;\dot{B}_{2,1}^{\frac {d}{2}} )}
+\big\| \q v^\ell\big\|_{\widetilde{L}^\infty(0,\infty;\dot{B}_{2,1}^{-1+\frac {d}{2}} )}+\|\nu a^h\|_{\widetilde{L}^\infty(0,\infty;\dot B^{\frac  dp}_{p,1})}\nonumber\\
&\quad\quad+\|\q v^h\|_{ \widetilde{L}^\infty(0,\infty;\dot B^{-1+\frac  dp}_{p,1})}+\|a^h\|_{ L^1(0,\infty;\dot B^{\frac  dp}_{p,1})}
+\nu\big\|  a^\ell\big\|_{L^1(0,\infty;\dot{B}_{2,1}^{1+\frac {d}{2}} )}
\nonumber\\
&\quad\quad+\nu^2\big\|  a^\ell\big\|_{L^1(0,\infty;\dot{B}_{2,1}^{2+\frac {d}{2}} )}
+\nu\big\| \q v^\ell\big\|_{L^1(0,\infty;\dot{B}_{2,1}^{1+\frac {d}{2}} )}+\nu \|\q v^h\|_{L^1(0,\infty;\dot B^{1+\frac  dp}_{p,1})}\nonumber\\
    &\quad\leq   C\exp\big(C(M+M^2)\big)\Big(\big\|a_0^\ell\big\|_{\dot{B}_{2,1}^{-1+\frac {d}{2}} }+\nu\big\|  a_0^\ell\big\|_{\dot{B}_{2,1}^{\frac {d}{2}} }+\nu\| a_0^h\|_{\dot B^{\frac  dp}_{p,1}}\nonumber\\
    &\hspace{5.6cm}+\big\| \q v_0^\ell\big\|_{\dot{B}_{2,1}^{-1+\frac {d}{2}} }+\|\q v_0^h\|_{\dot B^{-1+\frac  dp}_{p,1}}+M^2+\mu^2 \Big).
\end{align*}
If $a_0=0$, then the convergence $(\rho,v)\to(1,V)$ is obtained in the following sense:
\begin{align*}
&\sqrt{\mu^{-1}\nu}\|\rho-1\|_{L^\infty(0,\infty;\dot B^{\frac  dp}_{p,1})}
+  \|\p v -V\|_{L^\infty(0,\infty;\dot B^{-1+\frac dp}_{p,1})}    \nonumber\\
  &\quad\quad+\mu\| \p v -V\|_{L^1(0,\infty;\dot B^{1+\frac  dp}_{p,1})}  +\|\p v_t -V_t\|_{L^1(0,\infty;\dot B^{-1+\frac  dp}_{p,1})}  \nonumber\\
  &\quad \leq C\sqrt{\mu\nu^{-1}}.
\end{align*}
\end{theorem}

\begin{remark} It should be noted that the zero index Besov space $\dot B^{0}_{2,1}(\R^2)$ in Theorem \ref{yiyang}
 has now been  extended  by the  Besov space  $\dot B^{-1+\frac {d}{p}}_{p,1}(\R^d)$, which may has  the negative   index  $-1+\frac  dp<0$.   This implies the global well-posedness of compressible Navier-Stokes equations with highly oscillatory initial velocity field $v_0$, of which a typical example (see \cite[Proposition 2.9]{chenqionglei}) is
\begin{align*}
v_0(x)=\sin\bigl(\frac {x_1} {\varepsilon}\bigr)\phi(x), \quad \phi(x)\in \mathcal{S}(\R^d),\quad  p>d \mbox{ and } \varepsilon>0.
\end{align*}
This function is subject to the estimate:
\begin{align*}
\|v_0^\ell\|_{\dot{B}_{2,1}^{-1+\frac {d}{2}}}+ \|v_0^h\|_{\dot B^{-1+\frac  dp}_{p,1}}\le C\varepsilon^{1-\frac dp},
\end{align*}
for $C$ a constant independent of $\varepsilon>0$.
\end{remark}

\begin{remark}
 If $p=d=2$, Theorem  \ref{dingli} is identical to    \cite[Theorem 1.1]{danchin2018}. Especially,
when $d=2$ and $2\le p <4$, according to  \cite[Proposition 3.1]{huangjingchi}, the quantity  $M$ in Theorem  \ref{dingli} can be expressed precisely as
$$M=C\|\p v_0\|_{\dot{B}_{p,1}^{-1+\frac 2p}}\Big(1+\|\p v_0\|_{\dot{B}_{p,1}^{-1+\frac 2p}}\Big)
\exp\Big(\frac {C}{\mu^2}\|\p v_0\|^2_{\dot{B}_{p,1}^{-1+\frac 2p}}\Big).$$
If $d\geq3,$  we can construct some  examples of large initial data for
\eqref{classns} generating global smooth solutions. One can refer for instance to \cite{chemin,ponce,xuhuan} and citations therein.
\end{remark}

\begin{remark}
Recently, Danchin and Mucha  \cite{danchin1710}  derived  the large volume viscosity limit to the inhomogeneous incompressible Navier-Stokes equations
from \eqref{m}
in the two dimensional torus $\mathcal{T}^2.$ In particular, they can handle large variations of density.
\end{remark}

\

\subsection{Decomposition of (\ref{m}) by the Leray projection}
Without loss of generality, we  fix  the shear viscosity $\mu=1$ throughout the paper.

  Theorem \ref{dingli} is based on a decomposition of (\ref{m}) by using  (\ref{classns}). Employ the Leray projection to decompose  the velocity solution  into the compressible part $\q v$ and the incompressible part  $\p v+V$ as
$$v=\q u+ \p u+V\ \ \ \mbox{ for }\ \ \ u:=v-V $$
with $V$  the global solution of \eqref{classns}.
A simple computation implies
\begin{align}\label{biaoji}
 \q u=\q v,  \quad\diverg \q u =\diverg u.
\end{align}

For  $\rho=1+a$, we rewrite the second equation of \eqref{m} as
\begin{align}\label{yuanfangcheng}
v_t+(1+a)(v\cdot\nabla v)-\Delta v-(\lambda+1)\nabla\diverg v+P'\nabla a=-av_t.
\end{align}
Applying $\q$ to (\ref{yuanfangcheng}) and  using \eqref{biaoji} and the assumption $P'(1)=1$, we  get compressible part of system \eqref{m}:
\begin{eqnarray}\label{c}
({Compressible} \ {part})\quad\left\{\begin{aligned}
&a_t+\diverg\q u  =-\diverg(a(u+V)),\\
&(\q u)_t-\nu\Delta\q u+\nabla a=-\q H_1,\\
&a|_{t=0}=a_0,\quad \q u|_{t=0}=\q v_0,
\end{aligned}\right.
\end{eqnarray}
with
\begin{align*}%
H_1:=&\  \underbrace{a\left({V} _t+\p u_t+(\q u_t+\nabla a)\right)}_{H_1^{(1)}}
+\underbrace{(1+a)(u+{V} )\cdot\nabla (u+{V})}_{H_1^{(2)}}  +\underbrace{(k(a)-a)\nabla a}_{H_1^{(3)}},
\end{align*}
and
$$ k(a)=P'(1+a)-P'(1)=P'(1+a)-1.$$

Applying $\p$ to the  first equation of \eqref{classns} and   \eqref{yuanfangcheng}, respectively, and then taking the  difference between the two resultant  equations, we have
\begin{align}\label{i}
(\p u)_t-\Delta\p u=-\p H_2,\quad
\p u|_{t=0}=0,
\end{align}
with
\begin{align*}%
H_2:=&\underbrace{a\left({V} _t+\p u_t+(\q u_t+\nabla a)\right)}_{H_2^{(1)}}+\underbrace{(1+a)\p u\cdot\nabla({V} +\q u)}_{H_2^{(2)}}+\underbrace{a(u+{V} )\cdot\nabla \p u}_{H_2^{(3)}}
\\
&+\underbrace{(u+{V} )\cdot\nabla \p u}_{H_2^{(4)}}+\underbrace{(1+a)({V} \cdot\nabla\q u+\q u\cdot\nabla {V} )}_{H_2^{(5)}}
+\underbrace{a (\q u\cdot\nabla \q u+ {V} \cdot\nabla {V} )}_{H_2^{(6)}}.
\end{align*}

\subsection{Scheme for the proof of Theorem \ref{dingli}}
The proof of  Theorem \ref{dingli} lies on the observation  that the flow becomes incompressible if the  volume viscosity $\lambda$ is  sufficiently large. This mechanism  comes from  strong dissipation on the potential part
of the velocity when $\lambda$ is large.

In the proof of the main result, we  separate the original system \eqref{m} into incompressible part and compressible part. As the incompressible part satisfies a freedom heat equations, the estimates are standard. However,
 much effort has to be made in the examination of the compressible part.
By careful analysis of the linear system of the compressible part, we can find  the density $a$  has different  smoothing effect in low frequency and high frequency parts.
 Compared to the $L^2$ case obtained by Danchin and Mucha  \cite{danchin2018},
 the $L^p$ setting gives rise to an extra  difficulty to be overcome. In fact, we cannot get the smoothing effect of density in the low frequency and the damping effect of density in the high frequency  by employing  the traditional  energy argument, which  relies   heavily on  a cancelation, because the cancelation   is only valid in the  $L^2$ framework.  To get rid of  this difficulty,
 we shall derive the smoothing effect of  density in the $L^2$ setting in low frequency and  the damping effect of density in the $L^p$ setting in the high frequency, respectively. In the low frequency, one can follow the method used in \cite{danchin2018} to derive the desired estimates, whereas in the high frequency,
we follow an elementary energy approach in terms of \textit{effective velocity} developed by
Haspot \cite{haspot}, using Hoff's viscous effective flux in \cite{hoff2006}.

 As the solutions of incompressible part of \eqref{m} are constructed in the $L^p$ setting, all the estimates involved in incompressible part must work in the $L^p$  framework, even for us to derive the smoothing effect of density  in the low frequency.
For example, the nonlinear  term  $ v\cdot\nabla v $  is decomposed  into the following four parts:
\begin{equation}\label{fed}
v\cdot\nabla v=\p v\cdot\nabla \p v+\p v\cdot\nabla \q v+\q v\cdot\nabla \p v+\q v\cdot\nabla \q v.
\end{equation}
As the incompressible part $\p v$  lies in an $L^p$-type space yet $v\cdot\nabla v$ is estimated in the $L^2$ setting.  Hence,  we need some  product laws and commutators estimates (see   Lemma \ref{xinqing}, \eqref{pengyou} and \eqref{pengyou2} ) in the Besov spaces to deal with the  terms involved in $\p v$ in \eqref{fed}. A similar difficulty arises for  the  compressible part in the high frequency, for example the term $\q v^h\cdot\nabla \q v^h$.

The global solution to be obtained in Theorem \ref{dingli} will be extended from the
local solution given in the following:
\begin{theorem}(\cite[Theorem 2.2]{danchin2014}) \label{dingli6}
Let  $2\le d$, $1<p<2d$  and $v_0\in\dot B^{-1+\frac {d}{p}}_{p,1}(\R^d)$.
Assume $a_0=\rho_0-1\in\dot B^{\frac {d}{p}}_{p,1}(\R^d)$ and $
\inf_{x}\rho_0(x)>0.$ Then there exists a  maximal time $T^\ast>0$ so that \eqref{m} has a unique solution $(a,v)$ on $[0, T^\ast )$ satisfying, for any $T\in (0,T^\ast) $,
\begin{align*}
&a\in C([0,T];\dot{B}_{p,1}^{{\frac dp}}(\R^d))\cap
\widetilde{L}^\infty_T(\dot{B}_{p,1}^{{\frac dp}}(\R^d)),
\\
&v\in C([0,T];\dot{B}_{p,1}^{{-1+\frac dp}}(\R^d))\cap
\widetilde{L}^\infty_T(\dot{B}_{p,1}^{{-1+\frac dp}}(\R^d))\cap
L^1_T(\dot{B}_{p,1}^{{\frac dp+1}}(\R^d)).
\end{align*}
Moreover, continuation  beyond $T^\ast$ is possible if
\begin{equation*}
\int_0^{T^\ast}\|\nabla v\|_{L^\infty}\,dt<\infty,\quad \|a\|_{L^\infty(0,T^\ast;\dot B^{\frac dp}_{p,1})}<\infty\ \hbox{ and }\
\inf_{(t,x)\in[0,T^\ast)\times\R^d} \rho(t,x)>0.
\end{equation*}
\end{theorem}

  Thus the existence of the global solution in  Theorem \ref{dingli} becomes to show the maximal time $T^\ast$ to be $\infty$. This requires   global-in-time a priori estimates. Based on the decomposition in the previous subsection, the estimation will be divided  into three subsections in Section \ref{section3} through the  estimate of  the incompressible part in the critical $L^p$ framework, the estimate of  the high frequency of the compressible part of \eqref{c} and the estimate of  the low frequency of the compressible part of \eqref{c}.

\section{ Littlewood-Paley theory and preliminaries }

In this section, we recall some basic facts on Littlewood-Paley theory (see \cite{bcd} for instance).
Let $\chi$ and $\varphi$ be two smooth radial functions valued in the interval $[0,1]$ so that
the support of $\chi$ is the ball $\{\xi\in\R^d: |\xi|\leq\frac {4}{3}\}$, the support of $\varphi$ is  the annulus $\{\xi\in\R^d: \frac {3}{4}\leq|\xi|\leq\frac {8}{3}\}$ and
 \begin{equation*}
 \sum_{j \in \mathbb{Z}} \varphi(2^{-j}\xi)=1, \hspace{0.5cm}\forall \xi \neq0.
 \end{equation*}
   Let $\mathcal F$ be the Fourier transform.  
The homogeneous dyadic blocks $\dot{\Delta}_{j}$ and the homogeneous low-frequency cutoff operators
$\dot{S}_{j}$ are defined for all $j\in\mathbb{Z}$ by
\begin{align*}
\dot{\Delta}_{j}u={\mathcal F}^{-1}(\varphi(2^{-j}\cdot) {\mathcal F}u),\quad\quad
\dot{S}_{j}u={\mathcal F}^{-1}(\chi(2^{-j}\cdot) {\mathcal F}u) .
\end{align*}
Denote by $\mathscr{S}_h^{'}(\R^d)$ the space of tempered distributions subject to the condition
$$
\lim_{j\rightarrow -\infty}\dot{S}_j u=0.
$$
Then we have the  decomposition
$$
u = \sum_{j\in \Z} \dot \Delta_j u \ \ \  \forall u\in\mathscr{S}_h^{'}(\mathbb{R}^d).
$$
\begin{definition}
Let $s\in\R$ and $1\leq p,r\leq\infty$. The homogeneous Besov space $\dot{B}_{p,r}^{s}(\R^d)$
consists of all the distributions $u\in \mathscr{S}_{h}^{'}(\R^d)$ such that
$$
\|u\|_{\dot{B}_{p,r}^{s}}:=\left\|\big(2^{js}\|\dot{\Delta}_{j}u\|_{L^{p}}\big)_{j\in\mathbb{Z}}\right\|_{\ell^{r}}
<\infty.
$$
\end{definition}
This definition implies the following properties (see \cite[p71]{bcd}).
\begin{lemma}\label{equivalent defn}
Let $s\in\mathbb{R}$, $1\leq p,r\leq\infty$ and $u\in\mathscr{S}_h^{'}(\mathbb{R}^d)$. Then  $u$ belongs to $\dot{B}_{p,r}^{s}(\mathbb{R}^d)$ if and only if there exists a sequence  $\{c_{j,r}\}_{j\in\mathbb{Z}}$ with  $c_{j,r}\geq0$ and $\|c_{j,r}\|_{\ell^{r}}=1$ such that
$$
\|\dot{\Delta}_{j}u\|_{L^p}\leq Cc_{j,r}2^{-js}\|u\|_{\dot{B}_{p,r}^{s}},
$$
for a constant $C>0$.
If $r=1$, we set $d_j=c_{j,1}$.
\end{lemma}
\begin{lemma}\label{qianrudingli} For $s\in \R$ and  $1\le p,r\le \infty$, then there hold the estimates
$$\|u\|_{\dot{B}_{p,r}^{s}}\lesssim  \|\nabla u\|_{\dot{B}_{p,r}^{s-1}},\ \ \ \ \
\|\nabla u\|_{\dot{B}_{p,r}^{s-1}}\lesssim \|u\|_{\dot{B}_{p,r}^{s}}.$$
Here and in what follows, $a\lesssim b$ means the inequality  $a \le C b$ for a generic constant $C$.

Moreover, if $1\le p_1<p_2\le \infty$ and $1\le r_1<r_2\le \infty$, then we have
$${\dot{B}_{p_1,r_1}^{s}}(\R^d)\hookrightarrow {\dot{B}_{p_2,r_2}^{s-\frac{d}{p_1}-\frac{d}{p_2}}}(\R^d).$$
\end{lemma}

\begin{definition}
Let $s\in\mathbb{R}$ and $0<T\leq \infty$. The norm of the  Chemin-Lerner type Besov space is defined as
$$
\|u\|_{\widetilde{L}_{T}^{q}(\dot{B}_{p,1}^{s})}:=
\sum_{j\in\mathbb{Z}}
2^{jrs}\|\dot{\Delta}_{j}u\|_{L^q(0,T; L^p(\R^d))}
$$
for $1\leq p \leq \infty$ and $ 1\leq q <\infty$.
\end{definition}
This definition implies the inequality
\begin{align*}
\|u\|_{L^q_{T}(\dot{B}_{p,1}^s)}\le\|u\|_{\widetilde{L}^q_{T}(\dot{B}_{p,1}^s)},\hspace{0.5cm} \mathrm{for }\hspace{0.2cm}  q,\,p\ge 1.
\end{align*}
and the following interpolation property.
\begin{lemma}\label{neicha} (see \cite{bcd})
For $0<s_1<s_2,$ $0\leq \theta\leq 1$ and  $1\leq p,\, q_1,\, q_2\leq\infty$, then we have
\begin{align*}
\|u\|_{\widetilde{L}^q_{T}(\dot{B}_{p,1}^s)}\le\|u\|^\theta_{\widetilde{L}^{q_1}_{T}(\dot{B}_{p,1}^{s_1})}
\|u\|^{1-\theta}_{\widetilde{L}^{q_2}_{T}(\dot{B}_{p,1}^{s_2})} \ \  \mbox{ with }\ \ \frac {1}{q}=\frac {\theta}{q_1}+\frac {1-\theta}{q_2}, \ \ \ s=\theta s_1+(1-\theta)s_2.
\end{align*}
\end{lemma}

\begin{lemma}\label{bernstein}(Bernstein inequalities \cite{bcd})
Let $\mathcal{B}$ be a ball and $\mathcal{C}$ an annulus of $\mathbb{R}^d$ centered at the origin. For an integer $0\leq k\leq 2$ and reals
$1\le p \le q\le\infty$, there hold
\begin{align*}
&\sup_{|\alpha|=k}\|\partial^{\alpha}u\|_{L^q}\lesssim \sigma ^{k+d(\frac1p-\frac1q)}\|u\|_{L^p}, \ \mbox{ if } \ \ \mathrm{Supp} \,{\mathcal F}{u}\subset\sigma  \mathcal{B},\\
&\sigma ^k\|u\|_{L^p}\lesssim \sup_{|\alpha|=k}\|\partial^{\alpha}u\|_{L^p}
\lesssim \sigma ^{k}\|u\|_{L^p}, \ \mbox{ if } \ \ \mathrm{Supp} \,{\mathcal F}{u}\subset\sigma  \mathcal{C}
 \end{align*} with respect to scaling parameter $\sigma>0$.
\end{lemma}

The Bony decomposition  is
very effective in the estimate of nonlinear terms in fluid motion equations.
Here, we recall the decomposition in the homogeneous context:
\begin{align*}
uv=\dot{T}_uv+\dot{T}_vu+\dot{R}(u,v),
\end{align*}
where
$$\dot{T}_uv:=\sum_{j\in \Z}\dot{S}_{j-1}u\dot{\Delta}_jv \ \ \mbox{ and }\ \ \hspace{0.5cm}\dot{R}(u,v):=\sum_{j\in \Z}
\dot{\Delta}_ju\widetilde{\dot{\Delta}}_jv \  \ \mbox{ with } \ \ \widetilde{\dot{\Delta}}_jv:=\sum_{|j-j'|\le1}\dot{\Delta}_{j'}v.$$

The estimates of nonlinear terms of the compressible and incompressible equations are essentially based on the following lemmas.
\begin{lemma}\label{fenjie}(see \cite{bcd})
Let $s, s_1, s_2 \in\R$,  $1\leq p_1,p_2, r_1, r_2\leq\infty$,   $\frac  1p=\frac  1{p_1}+\frac 1{p_2}$, $ \frac1r =\frac1{r_1}+\frac1{r_2}$ and $\tau <0$.
Then we have
\begin{align*}
&\|\dot{T}_uv\|_{\dot B^s_{p,r}}\lesssim \|u\|_{L^{p_1}}\|v\|_{\dot B^{s}_{p_2,r}},
\quad\|\dot{T}_uv\|_{\dot B^{s+\tau}_{p,r}}\lesssim \|u\|_{\dot B^\tau_{p_1,\infty}}\|v\|_{\dot B^{s}_{p_2,r}},
\\
 &\|\dot{R}(u,v)\|_{\dot B^{s_1+s_2}_{p,r}}\lesssim\|u\|_{\dot B^{s_1}_{p_1,r_1}}
\|v\|_{\dot B^{s_2}_{p_2,r_2}}\,\,\mbox{ for }\,\, s_1+s_2>0,
\\
 &\|\dot{R}(u,v)\|_{\dot B^{0}_{p,\infty}}\lesssim\|u\|_{\dot B^{s_1}_{p_1,r_1}}
\|v\|_{\dot B^{s_2}_{p_2,r_2}}\,\,\mbox{ for }\,\, s_1+s_2=0.
\end{align*}
\end{lemma}

\begin{lemma}\label{daishu}(see \cite[Proposition A.1]{xujiang})
Let $d\ge 2$,  $1\leq p, q\leq \infty$, $s_1\leq \frac {d}{q}$, $s_2\leq d\min\{\frac 1p,\frac 1q\}$ and $s_1+s_2>d\max\{0,\frac 1p +\frac 1q -1\}$. Then we have,  for $ (u,v)\in\dot{B}_{q,1}^{s_1}({\mathbb R} ^d)\times\dot{B}_{p,1}^{s_2}({\mathbb R} ^d)$,
\begin{align*}
\|uv\|_{\dot{B}_{p,1}^{s_1+s_2 -\frac {d}{q}}}\lesssim \|u\|_{\dot{B}_{q,1}^{s_1}}\|v\|_{\dot{B}_{p,1}^{s_2}}.
\end{align*}
\end{lemma}

\begin{lemma}(see   \cite[Lemma 2.100]{bcd})\label{jiaohuanzi}
Let $d\ge 2$, $1\leq p, q\leq \infty$, $s\leq 1+d\min\{\frac 1p,\frac 1q\}$,
$v \in \dot{B}_{q,1}^{s}(\R^d)$ and $u\in \dot{B}_{p,1}^{\frac {d}{p}+1}(\R^d)$.
Assume that
$$s>-d\min\left\{\frac 1p,1-\frac {1}{q}\right\},\ \ or\ \ s>-1-d\min\left\{\frac 1p,1-\frac {1}{q}\right\}\ \ if\ \ \diverg u=0.$$
Then there holds the commutator estimate
$$
\|[u\cdot \nabla ,\dot{\Delta}_j]v \|_{L^q}\lesssim d_j 2^{-js}\|u\|_{\dot{B}_{p,1}^{\frac {d}{p}+1}}\|v \|_{\dot{B}_{q,1}^{s}},
$$
where and in what follows, we use the commutator symbol $[A,B]=AB-BA$ of  operators $A$ and $B$. 
\end{lemma}

The following estimates are implied from \cite[Lemma 2.16]{zhaixiaoping}.
\begin{lemma}\label{xinqing}
Let  $
2\leq p \leq \min\{4,2d/(d-2)\}$ for $ d>2$ and $2\leq p <4$ for $d=2$. Assume $A(D)$  a zero-order Fourier multiplier. For   $v^\ell\in\dot{B}^{-1+\frac{d}{2}}_{2,1}(\R^d)$,
$v^h\in\dot{B}^{-1+\frac  dp}_{p,1}(\R^d)$ and $\nabla u\in \dot{B}^{\frac dp}_{p,1}(\R^d)$, we have
\begin{align*}
&\sum_{j\le j_0}2^{(-1+\frac{d}{2})j}
\big\|\ddj([ A(D), u\cdot\nabla]v)\big\|_{L^2}\leq C (\big\|\nabla u^\ell\big\|_{\dot{B}^{\frac{d}{2}}_{2,1}}+\big\|\nabla u^h\big\|_{\dot{B}^{\frac{d}{p}}_{p,1}})(\big\|v^\ell\big\|_{\dot{B}^{ -1+\frac d2}_{2,1}}+\big\|v^h\big\|_{\dot{B}^{-1+\frac dp}_{p,1}}), \\
&\sum_{j\le j_0}2^{(-1+\frac{d}{2})j}
\big\|\ddj([ A(D), u\cdot\nabla]v)\big\|_{L^2}\leq C \big\|\nabla u\big\|_{\dot{B}^{\frac dp}_{p,1}}(\big\|v^\ell\big\|_{\dot{B}^{-1+\frac{d}{2}}_{2,1}}+\big\|v^h\big\|_{\dot{B}^{-1+\frac{d}{p}}_{p,1}}), \quad \hbox{if\  $\diverg u=0$},
\end{align*}
for a constant dependent on $j_0$.
\end{lemma}

The following lemma will be  used to get appropriate estimates of the solutions.
\begin{lemma}\label{dahai}  (see   \cite[Lemma 6.1]{helingbing})
Let $A(D)$ be a zero-order Fourier multiplier.
Let $j_0\in\Z,$  $\tau \in\R$,  $1\leq p_1,\ p_2\leq\infty$ and  $\frac  1p=\frac  1{p_1}+\frac  1{p_2}.$
Then we have 
\begin{align*}
\big\|[\dot S_{j_0} A(D), T_f ]g\big\|_{\dot B^{\tau +s}_{p,1}}
\leq C\|\nabla f \|_{\dot B^{s-1}_{p_1,1} }\|g\|_{\dot B^\tau _{p_2,\infty}},\quad\hbox{$s<1$,}
\\
\big\|[\dot S_{j_0} A(D), T_f ]g\big\|_{\dot B^{\tau +1}_{p,1}}
\leq C\|\nabla f \|_{L^{p_1}}\|g\|_{\dot B^\tau _{p_2,1}},\quad\hbox{s=1,}
\end{align*}
for a constant $C$ dependent on $j_0$.
\end{lemma}

Finally, we recall a composition result and a heat flow optimal regularity estimate.
  \begin{proposition}\label{fuhe2} (see \cite{bcd})
   Let $G$ with $G(0)=0$ be a smooth function defined on an open interval $I$
of $\R$ containing~$0.$
Then  the following estimates
$$
\|G(f)\|_{\dot B^{s}_{p,1}}\lesssim\|f\|_{\dot B^s_{p,1}}\quad\hbox{and}\quad
\|G(f)\|_{\widetilde{L}^q_T(\dot B^{s}_{p,1})}\lesssim\|f\|_{\widetilde{L}^q_T(\dot B^s_{p,1})}
$$
hold true for  $s>0,$ $1\leq p,\, q\leq \infty$ and
 $f$  valued in a bounded interval $J\subset I.$
\end{proposition}

\begin{proposition} (see \cite{bcd})\label{heat}
Let $\tau \in \R$, $\mu>0$, $T>0$, $1\leq p\leq\infty$ and $1\leq q_{2}\leq q_{1}\leq\infty$.  Let $u$  satisfy the heat equation
$$\left\{\begin{array}{lll}\partial_tu-\mu\Delta u=f,\\
u_{|t=0}=u_0.\end{array}
\right.$$
Then  the following a priori estimate
\begin{align*}
\mu^{\frac 1{q_1}}\|u\|_{\widetilde L_{T}^{q_1}(\dot B^{\tau +\frac 2{q_1}}_{p,1})}\lesssim
\|u_0\|_{\dot B^\tau _{p,1}}+\mu^{\frac {1}{q_2}-1}\|f\|_{\widetilde L^{q_2}_{T}(\dot B^{\tau -2+\frac 2{q_2}}_{p,1})}
\end{align*}holds true.
\end{proposition}

\section{Proof of Theorem \ref{dingli}}\label{section3}
The proof is to be completed through three  subsections with respect to  the incompressible part, the compressible part and their combination.

Let us begin with the  notation:
\begin{align*}%
\x:=&\big\|(a^\ell,\nu \nabla a^\ell, \q u^\ell)\big\|_{L^\infty(0,T;\dot{B}_{2,1}^{-1+\frac {d}{2}} )}
+\big\|\nu a^h\big\|_{L^\infty(0,T;\dot{B}_{p,1}^{\frac {d}{p}} )}+\big\| \q u^h\big\|_{L^\infty(0,T;\dot{B}_{p,1}^{-1+\frac {d}{p}} )},\nonumber\\
 \y:=&\big\|(\nu a^\ell,\nu^2 \nabla a^\ell, \nu\q u^\ell)\big\|_{L^1(0,T;\dot{B}_{2,1}^{1+\frac {d}{2}} )}+\big\| a ^h\big\|_{L^1(0,T;\dot{B}_{p,1}^{\frac {d}{p}} )}+\big\|\nu\q u^h\big\|_{L^1(0,T;\dot{B}_{p,1}^{1+\frac {d}{p}} )}\nonumber\\
 &+\big\|(\q u_t+\nabla a)^\ell\big\|_{L^1(0,T;\dot{B}_{2,1}^{-1+\frac {d}{2}} )}+\big\|(\q u_t+\nabla
 a)^h\big\|_{L^1(0,T;\dot{B}_{p,1}^{-1+\frac {d}{p}} )},\nonumber\\
\z:=&\big\|\p u\big\|_{L^\infty(0,T;\dot{B}_{p,1}^{-1+\frac {d}{p}} )}, \nonumber\\
\w:=&\big\|\p u_t\big\|_{L^1(0,T;\dot{B}_{p,1}^{-1+\frac {d}{p}} )}+\big\|\p u\big\|_{L^1(0,T;\dot{B}_{p,1}^{1+\frac {d}{p}} )},\nonumber\\
\mathcal{V}:=&  \|{V} \|_{L^\infty(0,T;\dot B^{-1+\frac {d}{p}}_{p,1})}+\|{V} _t\|_{L^1(0,T;\dot B^{-1+\frac {d}{p}}_{p,1})}+\| {V} \|_{L^1(0,T;\dot B^{1+\frac {d}{p}}_{p,1})}\le M.
\end{align*}

As assumed in Theorem \ref{dingli}, the system \eqref{classns} has  a unique global solution for  $\p v_0\in\dot{B}_{p,1}^{1+\frac {d}{p}}(\R^d).$ Thus the bound $M$ is  known.

We claim that if $\nu$ is large enough then one may find some (large) ${\eta}$ and (small) $\delta$
so that for all $T<T^*,$  the following bounds are  valid:
\begin{equation}\label{T1}
 \x+\y \leq {\eta}\quad\hbox{and}\quad
 \z+\w \leq \delta.
\end{equation}

\subsection{Estimates for the incompressible part of \eqref{m}}
Applying $\ddj$ to both side of \eqref{i} gives
\begin{align*}
(\dot{\Delta}_j\p u)_t-\Delta\dot{\Delta}_j\p u=-\dot{\Delta}_j\p H_2.
\end{align*}
Taking $L^2$ inner product with $|\ddj \p u|^{p-2}\ddj \p u$ to the above equation, we have
\begin{align}\label{Q1-1}
\frac 1p\frac {d}{dt}\|\dot{\Delta}_j\p u\|_{L^p}^p+C_12^{2j}\|\dot{\Delta}_j\p u\|_{L^p}^p\lesssim\|\dot{\Delta}_j\p H_2\|_{L^p}\|\dot{\Delta}_j\p u\|_{L^p}^{p-1},
\end{align}
in which we have used the  following fact  \cite[Appendix]{danchin2001}:
$$
-\int_{\R^d}\Delta\dot{\Delta}_j\p u\cdot|\dot{\Delta}_j\p u|^{p-2}\dot{\Delta}_j\p udx\ge C_12^{2j}\|\dot{\Delta}_j\p u\|_{L^p}^p $$
for some positive constant $C_1>0$.
Integrating  from $0$ to $T$ and using the H\"older inequality, we  get from \eqref{Q1-1} that
\begin{align}\label{Q2}
&\|\p u\|_{\widetilde{L}^\infty_T(\dot{B}_{p,1}^{-1+\frac {d}{p}} )}+\|\p u\|_{L^1_T(\dot{B}_{p,1}^{1+\frac {d}{p}} )}
\lesssim\int_0^T\|\p H_2\|_{\dot{B}_{p,1}^{-1+\frac {d}{p}} }dt.
\end{align}

\noindent In the following, we will deal with each terms in $\p H_2$.

Firstly, by Lemma \ref{daishu}, we have
\begin{align}\label{Q3}
\int_0^T\big\|\p H_2^{(1)}\big\|_{\dot{B}_{p,1}^{-1+\frac {d}{p}} }dt&\lesssim\int_0^T\|a(V_t+\p u_t+\q u_t+\nabla a)\|_{\dot{B}_{p,1}^{-1+\frac {d}{p}} }dt\nonumber\\
&\lesssim\nu^{-1} \|\nu a \|_{L^\infty_T(\dot{B}_{p,1}^{\frac {d}{p}} )}\|(\q u_t+\nabla a,\p u_t,V_t)\|_{L^1_T(\dot{B}_{p,1}^{-1+\frac {d}{p}} )}\nonumber\\
&\lesssim\nu^{-1} \Big(\|\nu a^h \|_{L^\infty_T(\dot{B}_{p,1}^{\frac {d}{p}} )}+\|\nu a^\ell \|_{L^\infty_T(\dot{B}_{2,1}^{\frac {d}{2}} )}\Big)\Big(\|(\q u_t+\nabla a)^\ell\|_{L^1_T(\dot{B}_{2,1}^{-1+\frac {d}{2}} )}\nonumber\\
&\quad\quad\quad\quad+\|((\q u_t+\nabla a)^h,\p u_t,V_t)\|_{L^1_T(\dot{B}_{p,1}^{-1+\frac {d}{p}} )}\Big)\nonumber\\
&\lesssim \nu^{-1}\x(\y+\w+\v).
\end{align}
Similarly, we have
\begin{align}\label{Q5}
\int_0^T\big\|\p H_2^{(2)}\big\|_{\dot{B}_{p,1}^{-1+\frac {d}{p}} }dt
\lesssim&(1+\|a\|_{L^\infty_T(\dot{B}_{p,1}^{\frac {d}{p}}) })\int_0^T(\|\nabla V\|_{\dot{B}_{p,1}^{\frac {d}{p}} }+\|\nabla\q u\|_{\dot{B}_{p,1}^{\frac {d}{p}} })\|\p u\|_{\dot{B}_{p,1}^{-1+\frac {d}{p}} }dt,
\end{align}
and
\begin{align}\label{Q6}
&\int_0^T\big\|\p H_2^{(3)}\big\|_{\dot{B}_{p,1}^{-1+\frac {d}{p}} }dt
\lesssim\nu^{-1} \x(\z+\x+\v)\w.
\end{align}
By the interpolation inequality in Lemma \ref{daishu}, we get
\begin{align}\label{Q8}
\int_0^T\big\|\p H_2^{(4)}\big\|_{\dot{B}_{p,1}^{-1+\frac {d}{p}} }dt\lesssim&\int_0^T\|(u+{V} )\cdot\nabla \p u\|_{\dot{B}_{p,1}^{-1+\frac {d}{p}} }dt\nonumber\\
\lesssim&\int_0^T(\|u\|_{\dot{B}_{p,1}^{\frac {d}{p}} }+\|V\|_{\dot{B}_{p,1}^{\frac {d}{p}} })\|\nabla \p u\|_{\dot{B}_{p,1}^{-1+\frac {d}{p}} }dt\nonumber\\
\lesssim&\int_0^T(\|u\|_{\dot{B}_{p,1}^{\frac {d}{p}} }+\|V\|_{\dot{B}_{p,1}^{\frac {d}{p}} })\| \p u\|_{\dot{B}_{p,1}^{-1+\frac {d}{p}} }^{\frac 12}\| \p u\|_{\dot{B}_{p,1}^{1+\frac {d}{p}} }^{\frac 12}dt\nonumber\\
\lesssim&\varepsilon\| \p u\|_{L^1_T(\dot{B}_{p,1}^{1+\frac {d}{p}} )}+\int_0^T(\|u\|_{\dot{B}_{p,1}^{\frac {d}{p}} }^2+\|V\|_{\dot{B}_{p,1}^{\frac {d}{p}} }^2)\| \p u\|_{\dot{B}_{p,1}^{-1+\frac {d}{p}} }dt,
\end{align}
\begin{align}\label{Q9}
&\int_0^T\big\|\p H_2^{(5)}\big\|_{\dot{B}_{p,1}^{-1+\frac {d}{p}} }dt\nonumber\\
&\quad\lesssim\int_0^T(1+\|a\|_{\dot{B}_{p,1}^{\frac {d}{p}} })\|\q u\|_{\dot{B}_{p,1}^{-1+\frac {d}{p}} }^{\frac 12}\|\q u\|_{\dot{B}_{p,1}^{1+\frac {d}{p}} }^{\frac 12}\|V\|_{\dot{B}_{p,1}^{-1+\frac {d}{p}} }^{\frac 12}\|V\|_{\dot{B}_{p,1}^{1+\frac {d}{p}} }^{\frac 12}dt
\nonumber\\
&\quad\lesssim(1+\|a\|_{L^\infty_T(\dot{B}_{p,1}^{\frac {d}{p}} )})\|\q u\|_{L^\infty_T(\dot{B}_{p,1}^{-1+\frac {d}{p}} )}^{\frac 12}\|\q u\|_{L^1_T(\dot{B}_{p,1}^{1+\frac {d}{p}} )}^{\frac 12}\|V\|_{L^\infty_T(\dot{B}_{p,1}^{-1+\frac {d}{p}} )}^{\frac 12}\|V \|_{L^1_T(\dot{B}_{p,1}^{1+\frac {d}{p}} )}^{\frac 12}\nonumber\\
&\quad\lesssim(1+\nu^{-1} \x)\nu^{-\frac 12}\x^{\frac 12}\y^{\frac 12}\v,
\end{align}
\begin{align}\label{Q10}
&\int_0^T\big\|\p H_2^{(6)}\big\|_{\dot{B}_{p,1}^{-1+\frac {d}{p}} }dt\nonumber\\
&\quad\lesssim
\|a(\q u\cdot\nabla \q u+V\cdot\nabla V)\|_{L^1_T(\dot{B}_{p,1}^{-1+\frac {d}{p}} )}\nonumber\\
&\quad\lesssim\|a\|_{L^\infty_T(\dot{B}_{p,1}^{\frac {d}{p}} )}(\|\q u\|_{L^\infty_T(\dot{B}_{p,1}^{-1+\frac {d}{p}} )}\|\q u\|_{L^1_T(\dot{B}_{p,1}^{1+\frac {d}{p}} )}+\|V\|_{L^\infty_T(\dot{B}_{p,1}^{-1+\frac {d}{p}} )}\|V \|_{L^1_T(\dot{B}_{p,1}^{1+\frac {d}{p}} )})\nonumber\\
&\quad\lesssim\nu^{-1} \x(\nu^{-1}\x\y+\v^2).
\end{align}

Assuming from now on that
\begin{equation}\label{tiaojian1}
\nu^{-1}\eta\ll 1.
\end{equation}
After a simple computation,   it follows  from definition of $\x$ and \eqref{T1} that
\begin{align}\label{axiao}
\|a\|_{L^\infty_T(\dot{B}_{p,1}^{\frac {d}{p}} )}\lesssim \nu^{-1}\|\nu a^\ell\|_{L^\infty_T(\dot{B}_{2,1}^{\frac {d}{2}})}+\nu^{-1}\|\nu a^h\|_{L^\infty_T(\dot{B}_{p,1}^{\frac {d}{p}} )}
\lesssim \nu^{-1}\eta\ll1.
\end{align}
Inserting the estimates \eqref{Q3}--\eqref{Q10} into \eqref{Q2} and choosing $\varepsilon $ small enough, we have
\begin{align}\label{Q11}
&\|\p u\|_{\widetilde{L}^\infty_T(\dot{B}_{p,1}^{-1+\frac {d}{p}} )}+\|\p u\|_{L^1_T(\dot{B}_{p,1}^{1+\frac {d}{p}} )}\nonumber\\
&\quad\lesssim\int_0^T\Big(\|(\nabla V,\nabla \q u)\|_{\dot{B}_{p,1}^{\frac {d}{p}} }+\|u\|_{\dot{B}_{p,1}^{\frac {d}{p}} }^2+\|V\|_{\dot{B}_{p,1}^{\frac {d}{p}} }^2\Big)\|\p u\|_{\dot{B}_{p,1}^{-1+\frac {d}{p}} } dt+\nu^{-\frac 12}\x^{\frac 12}\y^{\frac 12}\v\nonumber\\
&\quad\quad
+\nu^{-1} \x(\z+\x+\v)\w+\nu^{-1}(\y+\w+\v)\x+\nu^{-1} \x(\v^2+\nu^{-1}\x\y).
\end{align}
By \eqref{i}, we have
\begin{align}\label{Q12-2}
\big\|\p u_t\big\|_{L^1_T(\dot{B}_{p,1}^{-1+\frac {d}{p}} )}
\lesssim\|\p u\|_{L^1_T(\dot{B}_{p,1}^{1+\frac {d}{p}} )}+\int_0^T\|\p H_2\|_{\dot{B}_{p,1}^{-1+\frac {d}{p}} }dt.
\end{align}
The combination of \eqref{Q11} and  \eqref{Q12-2} with the  Gronwall inequality  produces that
\begin{align}\label{Q12}
\z+\w
\le&\exp\Big(C\int_0^T\Big(\|(\nabla V,\nabla \q u)\|_{\dot{B}_{p,1}^{\frac {d}{p}} }+\|u\|_{\dot{B}_{p,1}^{\frac {d}{p}} }^2+\|V\|_{\dot{B}_{p,1}^{\frac {d}{p}} }^2\Big)dt\Big)\Big\{\nu^{-\frac 12}\x^{\frac 12}\y^{\frac 12}\v\nonumber\\
&\quad\quad\quad\quad\quad+\nu^{-1} \x(\z+\x+\v)\w+\nu^{-1}(\y+\w+\v)\x+\nu^{-1} \x\v^2\Big\}.
\end{align}

\bigskip

\subsection{High frequencies for the compressible part of \eqref{m}}
To estimate the high frequencies of $(a,\q u),$ we follow the approach of \cite{haspot} and  introduce
 the following ``effective'' velocity field
 $$w=\q u+\nu^{-1}(-\Delta)^{-1}\nabla a.$$
Multiplying by $\nu^{-1}(-\Delta)^{-1}\nabla$  on the first equation in \eqref{c} and then adding the resultant equation to the
  second one in \eqref{c}, we deduce that
\begin{align}\label{A1}
w_t-\nu \Delta w=&\nu^{-1}\q u+\nu^{-1}\q(a(u+V))-\q H_1\nonumber\\
=&\nu^{-1}w-\nu^{-2}(-\Delta)^{-1}\nabla a+\nu^{-1}\q(a(u+V))-\q H_1.
\end{align}

Applying the operator $\ddj$ on  \eqref{A1} and multiplying by $|\ddj w|^{p-2}\ddj w$ to the resultant equation,  we get that
\begin{align}\label{Q1-12322}
\frac 1p\frac {d}{dt}\|\dot{\Delta}_j w\|_{L^p}^p+C\nu2^{2j}\|\dot{\Delta}_j w\|_{L^p}^p
\lesssim&\nu^{-1}\|\dot{\Delta}_jw\|_{L^p}^p+\|\nu^{-2}(-\Delta)^{-1}\nabla a\|\dot{\Delta}_jw\|_{L^p}^{p-1}\nonumber\\
&+
\|(\nu^{-1}\q(a(u+V))-\q H_1)\|_{L^p}\|\dot{\Delta}_jw\|_{L^p}^{p-1}.
\end{align}
Hence multiplying  \eqref{Q1-12322} by  $2^{(-1+\frac dp)j}/\|\dot{\Delta}_jw\|_{L_p}^{p-1}$, then integrating  with respect to $t$  and summing up the resultant equations  for the high frequencies $\dot \Delta_j w$  only, we get
\begin{align}\label{A2}
\|w\|^h_{ \widetilde{L}_T^\infty(\dot B^{-1+\frac  dp}_{p,1})}+\nu \|w\|^h_{L^1_T(\dot B^{1+\frac  dp}_{p,1})}
\lesssim& \|w_0\|^h_{\dot B^{-1+\frac  dp}_{p,1}}
+ \|\nu^{-1}w\|^h_{L^1_T(\dot B^{-1+\frac  dp}_{p,1})}+ \|\nu^{-2}(-\Delta)^{-1}\nabla a\|^h_{L^1_T(\dot B^{-1+\frac  dp}_{p,1})}\nonumber\\
&+ \|\nu^{-1}\q(a(u+V))\|^h_{L^1_T(\dot B^{-1+\frac  dp}_{p,1})}+\|\q H_1\|^h_{L^1_T(\dot B^{-1+\frac  dp}_{p,1})}\nonumber\\
\lesssim&  \|w_0\|^h_{\dot B^{-1+\frac  dp}_{p,1}}
+ \nu^{-1}\|w\|^h_{L^1_T(\dot B^{-1+\frac  dp}_{p,1})}+\nu^{-2}2^{-2j_0}\|a^h\|_{L^1_T(\dot B^{\frac  dp}_{p,1})}\nonumber\\
&+ \|\nu^{-1}a(u+V)\|^h_{L^1_T(\dot B^{-1+\frac  dp}_{p,1})}+ \|\q H_1\|^h_{L^1_T(\dot B^{-1+\frac  dp}_{p,1})}.
\end{align}
By Lemma \ref{daishu} and Young's inequality, we have
\begin{align}\label{A3}
 \|\nu^{-1}a(u+V)\|^h_{L^1_T(\dot B^{-1+\frac  dp}_{p,1})}
 & \lesssim \nu^{-1}\|a\|_{L^2_T(\dot B^{\frac  dp}_{p,1})}\|u+V\|_{L^2_T(\dot B^{\frac  dp}_{p,1})}\nonumber\\
 & \lesssim\|a\|_{L^2_T(\dot B^{\frac  dp}_{p,1})}^2+\nu^{-2}\|(u,V)\|_{L^2_T(\dot B^{\frac  dp}_{p,1})}^2\nonumber\\
 &\lesssim\int_0^T\Big(\nu^{-1}\|a^h\|_{\dot B^{\frac  dp}_{p,1}}+\nu^{-1}\|\nu a^\ell\|_{\dot B^{1+\frac  d2}_{2,1}}\Big)\Big(\|\nu a^h\|_{\dot B^{\frac  dp}_{p,1}}+\| a^\ell\|_{\dot B^{-1+\frac  d2}_{2,1}}\Big)dt\nonumber\\
 &  \quad +\nu^{-2}\z\w+\nu^{-3}\x\y+\nu^{-2}\v^2,
\end{align}
where we have used the following estimate:
\begin{align}\label{A5}
&\|(u,V)\|_{L^2_T(\dot B^{\frac  dp}_{p,1})}^2\nonumber\\
&\quad\lesssim\int_0^T\Big( \|\p u\|_{\dot B^{\frac  dp}_{p,1}}^2+\|\q u^h\|_{\dot B^{\frac  dp}_{p,1}}^2+\|\q u^\ell\|_{\dot B^{\frac  d2}_{2,1}}^2+\|V\|_{\dot B^{\frac  dp}_{p,1}}^2 \Big)dt\nonumber\\
&\quad\lesssim\int_0^T\Big( \|\p u\|_{\dot B^{-1+\frac  dp}_{p,1}}\|\p u\|_{\dot B^{1+\frac  dp}_{p,1}}+\|\q u^h\|_{\dot B^{-1+\frac  dp}_{p,1}}\|\q u^h\|_{\dot B^{1+\frac  dp}_{p,1}}\nonumber\\
&\quad\quad+\|\q u^\ell\|_{\dot B^{-1+\frac  d2}_{2,1}}\|\q u^\ell\|_{\dot B^{1+\frac  d2}_{2,1}}+\|V\|_{\dot B^{-1+\frac  dp}_{p,1}}\|V\|_{\dot B^{1+\frac  dp}_{p,1}} \Big)dt\nonumber\\
&\quad\lesssim\|\p u\|_{ \widetilde{L}_T^\infty(\dot B^{-1+\frac  dp}_{p,1})} \|\p u\|_{L^1_T(\dot B^{1+\frac  dp}_{p,1})}+\nu^{-1}\|\q u^h\|_{ \widetilde{L}_T^\infty(\dot B^{-1+\frac  dp}_{p,1})} \|\nu\q u^h\|_{L^1_T(\dot B^{1+\frac  dp}_{p,1})}\nonumber\\
&\quad\quad+\nu^{-1}\|\q u^\ell\|_{ \widetilde{L}_T^\infty(\dot B^{-1+\frac  d2}_{2,1})} \|\nu\q u^\ell\|_{L^1_T(\dot B^{\frac  d2+1}_{2,1})}+\|V\|_{ \widetilde{L}_T^\infty(\dot B^{-1+\frac  dp}_{p,1})} \|V\|_{L^1_T(\dot B^{1+\frac  dp}_{p,1})}\nonumber\\
&\quad\lesssim\z\w+\nu^{-1}\x\y+\v^2.
\end{align}
Thanks to Lemma \ref{daishu}, we obtain that
\begin{align}\label{A6}
 \big\|\q H_1^{(1)}\big\|^h_{L^1_T(\dot B^{-1+\frac  dp}_{p,1})}
 &\lesssim \|a\left(\mathcal{V} _t+\p u_t+(\q u_t+\nabla a)\right)
\|^h_{L^1_T(\dot B^{-1+\frac  dp}_{p,1})}\nonumber\\
 &\lesssim\|a\|_{L^\infty_T(\dot B^{\frac  dp}_{p,1})}\|\left(\p u_t,\mathcal{V} _t,(\q u_t+\nabla a)\right)
\|_{L^1_T(\dot B^{-1+\frac  dp}_{p,1})}\nonumber\\
&\lesssim\nu^{-1}\x(\w+\v+\y).
\end{align}
Similarly, we have
\begin{align}\label{A8}
& \big\|\q H_1^{(2)}\big\|^h_{L^1_T(\dot B^{-1+\frac  dp}_{p,1})} \nonumber\\
 &\quad\lesssim(1+\|a\|_{L^\infty_T(\dot B^{\frac  dp}_{p,1})})\|u+\mathcal{V}\|_{L^\infty_T(\dot B^{-1+\frac  dp}_{p,1})}\|
\nabla (u+\mathcal{V}) \|_{L^1_T(\dot B^{\frac  dp}_{p,1})}\nonumber\\
&\quad\lesssim(1+\|a\|_{L^\infty_T(\dot B^{\frac  dp}_{p,1})})\Bigg\{\int_0^T\|(\nabla  u,\nabla V)\|_{\dot{B}_{p,1}^{\frac {d}{p}} }
\Big(\|\q u^h\|_{\dot{B}_{p,1}^{-1+\frac {d}{p}}}+\|\q u^\ell\|_{\dot{B}_{2,1}^{-1+\frac {d}{2}} }\Big)dt\nonumber\\
&\quad\quad+\int_0^T\Big(\|(\nabla  \p u,\nabla \q u^h,\nabla V)\|_{\dot{B}_{p,1}^{\frac {d}{p}} }+\|\nabla \q u^\ell\|_{\dot{B}_{2,1}^{\frac {d}{2}} }\Big)\|(\p u,V)\|_{\dot{B}_{p,1}^{-1+\frac {d}{p}})}dt\Bigg\}\nonumber\\
&\quad\lesssim(\v+\w)(\z+\v)+\nu^{-1}\y(\z+\v)\nonumber\\
&\quad\quad+\int_0^T\|(\nabla  u,\nabla V)\|_{\dot{B}_{p,1}^{\frac {d}{p}} }
\Big(\|\q u^h\|_{\dot{B}_{p,1}^{-1+\frac {d}{p}}}+\|\q u^\ell\|_{\dot{B}_{2,1}^{-1+\frac {d}{2}} }\Big)dt,
\end{align}
after the use of  the smallness assumption \eqref{axiao} on $\|a\|_{L^\infty_T(\dot B^{\frac  dp}_{p,1})}$.

From Proposition \ref{fuhe2}, we  get
\begin{align}\label{A9}
&\|(k(a)-a)\nabla a\|^h_{L^1_T(\dot B^{-1+\frac  dp}_{p,1})}\nonumber\\
&\quad\lesssim\int_0^T\Big(\nu^{-1}\|a^h\|_{\dot B^{\frac  dp}_{p,1}}+\nu^{-1}\|\nu a^\ell\|_{\dot B^{1+\frac  d2}_{2,1}}\Big)\Big(\|\nu a^h\|_{\dot B^{\frac  dp}_{p,1}}+\| a^\ell\|_{\dot B^{-1+\frac  d2}_{2,1}}\Big)dt.
\end{align}

Thus, the combination of  \eqref{A3}-- \eqref{A9} with  \eqref{A2} implies that
\begin{align}\label{A10}
&\|w\|^h_{ \widetilde{L}^\infty_T(\dot B^{-1+\frac  dp}_{p,1})}+\nu \|w\|^h_{L^1_T(\dot B^{1+\frac  dp}_{p,1})}
\nonumber\\
&\quad\lesssim \|w_0\|^h_{\dot B^{-1+\frac  dp}_{p,1}}
+ \nu^{-1}\|w\|^h_{L^1_T(\dot B^{-1+\frac  dp}_{p,1})}+\nu^{-2}\|a^h\|_{L^1_T(\dot B^{\frac  dp}_{p,1})}+(\v+\w)(\z+\v)\nonumber\\
&\quad\quad+\nu^{-1}\y(\z+\v)
+\nu^{-1}\x(\w+\v+\y)+\nu^{-2}\z\w+\nu^{-3}\x\y+\nu^{-2}\v^2\nonumber\\
&\quad\quad+\int_0^T\|(\nabla   u,\nabla V)\|_{\dot{B}_{p,1}^{\frac {d}{p}} }\Big(\|\q u^h\|_{\dot{B}_{p,1}^{-1+\frac {d}{p}})}+\|\q u^\ell\|_{\dot{B}_{2,1}^{-1+\frac {d}{2}} }\Big)dt\nonumber\\
&\quad\quad+\int_0^T\Big(\nu^{-1}\|a^h\|_{\dot B^{\frac  dp}_{p,1}}+\nu^{-1}\|\nu a^\ell\|_{\dot B^{1+\frac  d2}_{2,1}}\Big)\Big(\|\nu a^h\|_{\dot B^{\frac  dp}_{p,1}}+\| a^\ell\|_{\dot B^{-1+\frac  d2}_{2,1}}\Big)dt.
\end{align}


We now detect damping
phenomenon on the high frequencies of the  density.

As $\diverg\q u=\diverg w+\nu^{-1} a$, we deduce from the first equation in \eqref{c} that
\begin{align}
\partial_ta+(u+V)\cdot\nabla a+\nu^{-1}a=-a\,\diverg  u-\diverg  w.
\end{align}
To bound the high frequencies of $a,$ we write
\begin{eqnarray}\label{A11}
\partial_t\ddj a+(u+V)\cdot\nabla\ddj a+\nu^{-1}\ddj a
=[(u+V)\cdot\nabla,\ddj ]a-\ddj(a\,\diverg  u+\diverg  w).
\end{eqnarray}

\noindent
Taking $L^2$ inner product with $|\ddj a|^{p-2}\ddj a$, using integrating by part and  the H\"older inequality,  we thus get, for $t\geq0,$
\begin{align}\label{W1}
&\|\ddj a(t)\|_{L^p}+\nu^{-1}\int_0^T\|\ddj a\|_{L^p}\,dt\nonumber\\
&\quad\leq\|\ddj a_0\|_{L^p}
+\frac 1p\int_0^T\|\diverg (u+V)\|_{L^\infty}\|\ddj a\|_{L^p}\,dt
\nonumber\\
&\quad\quad+\int_0^T\|[(u+V)\cdot\nabla,\ddj ]a\|_{L^p}\,dt+\int_0^T\|\ddj(a\,\diverg  u+\diverg  w)\|_{L^p}\,dt
.
\end{align}
By Lemmas \ref{daishu} and \ref{jiaohuanzi}, we  get
\begin{align*}
&\|a\,\diverg  u\|_{\dot B^{\frac  dp}_{p,1}}\lesssim \|a\|_{\dot B^{\frac  dp}_{p,1}}\|\diverg  u\|_{\dot B^{\frac  dp}_{p,1}},\\
&\sum_{j\in\Z}2^{\frac  dpj}\|[(u+V)\cdot\nabla,\ddj ]a\|_{L^p}\leq C\|\nabla (u+V)\|_{\dot B^{\frac  dp}_{p,1}}\|a\|_{\dot B^{\frac  dp}_{p,1}}.
\end{align*}
Multiplying \eqref{W1} by $2^{\frac  dpj},$ using the  embedding relation ${\dot B^{\frac  dp}_{p,1}}(\R^d)\hookrightarrow L^\infty(\R^d)$, we have
\begin{align*}
\|a^h\|_{\widetilde{L}^\infty_T(\dot B^{\frac  dp}_{p,1})}+\nu^{-1}\|a^h\|_{ L^1_T(\dot B^{\frac  dp}_{p,1})}\le& \|a_0\|^h_{\dot B^{\frac  dp}_{p,1}}+\|w\|^h_{L^1_T(\dot B^{1+\frac  dp}_{p,1})}+C\int_0^T\|(\nabla u,\nabla V)\|_{\dot B^{\frac  dp}_{p,1}}\|a\|_{\dot B^{\frac  dp}_{p,1}}\,dt,
\end{align*}
which implies that
\begin{align}\label{W8}
\|\nu a\|^h_{\widetilde{L}^\infty_T(\dot B^{\frac  dp}_{p,1})}+\|a^h\|_{ L^1_T(\dot B^{\frac  dp}_{p,1})}
\le& \|\nu a_0\|^h_{\dot B^{\frac  dp}_{p,1}}+\nu\|w\|^h_{L^1_T(\dot B^{1+\frac  dp}_{p,1})}\nonumber\\
&+C\int_0^T\|(\nabla u,\nabla V)\|_{\dot B^{\frac  dp}_{p,1}}(\|\nu a^h\|_{\dot B^{\frac  dp}_{p,1}}+\|\nu a^\ell\|_{\dot B^{\frac  d2}_{2,1}})dt.
\end{align}
By \eqref{A10} and \eqref{W8} and  since  $\nu $ is large enough ( here we use $\nu>1$), we get
\begin{align}\label{W10}
&\|\nu a\|^h_{\widetilde{L}^\infty_T(\dot B^{\frac  dp}_{p,1})}+\|a^h\|_{ L^1_T(\dot B^{\frac  dp}_{p,1})}+\|w\|^h_{ \widetilde{L}^\infty_T(\dot B^{-1+\frac  dp}_{p,1})}+\nu \|w\|^h_{L^1_T(\dot B^{1+\frac  dp}_{p,1})}
\nonumber\\
&\quad\lesssim  \|\nu a_0\|^h_{\dot B^{\frac  dp}_{p,1}}+\|w_0\|^h_{\dot B^{-1+\frac  dp}_{p,1}}+(\v+\w)(\z+\v)+\nu^{-1}\y(\z+\v)\nonumber\\
&\quad\quad
+\nu^{-1}\x(\w+\v+\y)+\nu^{-2}\z\w+\nu^{-3}\x\y+\nu^{-2}\v^2\nonumber\\
&\quad\quad+\int_0^T\Big(\nu^{-1}\|a^h\|_{\dot B^{\frac  dp}_{p,1}}+\nu^{-1}\|\nu a^\ell\|_{\dot B^{1+\frac  d2}_{2,1}}\Big)\Big(\|\nu a^h\|_{\dot B^{\frac  dp}_{p,1}}+\| a^\ell\|_{\dot B^{-1+\frac  d2}_{2,1}}\Big)dt\nonumber\\
&\quad\quad
+\int_0^T\|(\nabla   u,\nabla V)\|_{\dot{B}_{p,1}^{\frac {d}{p}} }\Big(\|\q u^h\|_{\dot{B}_{p,1}^{-1+\frac {d}{p}}}+\|\q u^\ell\|_{\dot{B}_{2,1}^{-1+\frac {d}{2}} }+\|\nu a^h\|_{\dot B^{\frac  dp}_{p,1}}+\|\nu a^\ell\|_{\dot B^{\frac  d2}_{2,1}}\Big)dt.
\end{align}
Recalling  that
 $$w=\q u+\nu^{-1}(-\Delta)^{-1}\nabla a,$$
 we  get
\begin{align}\label{W11}
\|\q u\|^h_{ \widetilde{L}^\infty_T(\dot B^{-1+\frac  dp}_{p,1})}
\lesssim\|w\|^h_{ \widetilde{L}^\infty_T(\dot B^{-1+\frac  dp}_{p,1})}\!+\|\nu a\|^h_{ \widetilde{L}^\infty_T(\dot B^{\frac  dp}_{p,1})}\!,
\nu \|\q u\|^h_{L^1_T(\dot B^{1+\frac  dp}_{p,1})}
\lesssim\nu \|w\|^h_{L^1_T(\dot B^{1+\frac  dp}_{p,1})}\!+\|a^h\|_{ L^1_T(\dot B^{\frac  dp}_{p,1})}
\end{align}
and so,  by   \eqref{W10} and \eqref{W11},
\begin{align}\label{W12}
&\|\nu a\|^h_{\widetilde{L}^\infty_T(\dot B^{\frac  dp}_{p,1})}+\|a^h\|_{ L^1_T(\dot B^{\frac  dp}_{p,1})}+\|\q u\|^h_{ \widetilde{L}^\infty_T(\dot B^{-1+\frac  dp}_{p,1})}+\nu \|\q u\|^h_{L^1_T(\dot B^{1+\frac  dp}_{p,1})}
\nonumber\\
&\quad\lesssim  \|\nu a_0\|^h_{\dot B^{\frac  dp}_{p,1}}+\|w_0\|^h_{\dot B^{-1+\frac  dp}_{p,1}}+(\v+\w)(\z+\v)+\nu^{-1}\y(\z+\v)\nonumber\\
&\quad\quad
+\nu^{-1}\x(\w+\v+\y)+\nu^{-2}\z\w+\nu^{-3}\x\y+\nu^{-2}\v^2\nonumber\\
&\quad\quad+\int_0^T\Big(\nu^{-1}\|a^h\|_{\dot B^{\frac  dp}_{p,1}}+\nu^{-1}\|\nu a^\ell\|_{\dot B^{1+\frac  d2}_{2,1}}\Big)\Big(\|\nu a^h\|_{\dot B^{\frac  dp}_{p,1}}+\| a^\ell\|_{\dot B^{-1+\frac  d2}_{2,1}}\Big)dt\nonumber\\
&\quad\quad
+\int_0^T\|(\nabla   u,\nabla V)\|_{\dot{B}_{p,1}^{\frac {d}{p}} }\Big(\|\q u^h\|_{\dot{B}_{p,1}^{-1+\frac {d}{p}}}+\|\q u^\ell\|_{\dot{B}_{2,1}^{-1+\frac {d}{2}} }+\|\nu a^h\|_{\dot B^{\frac  dp}_{p,1}}+\|\nu a^\ell\|_{\dot B^{\frac  d2}_{2,1}}\Big)dt.
\end{align}
We employ  the second equation of \eqref{c} to produce that
\begin{align}\label{W121}
\big\|(\q u_t+\nabla
 a)^h\big\|_{L^1_T(\dot{B}_{p,1}^{-1+\frac {d}{p}} )}
\lesssim&\big\|\nu\q u^h\big\|_{L^1_T(\dot{B}_{p,1}^{1+\frac {d}{p}} )}+\big\|\q H_1\big\|_{L^1_T(\dot{B}_{p,1}^{-1+\frac {d}{p}} )}.
\end{align}
Thus, the combination of  \eqref{W12} and  \eqref{W121} implies  that
\begin{align}\label{W123}
&\|\nu a\|^h_{\widetilde{L}^\infty_T(\dot B^{\frac  dp}_{p,1})}+\|a^h\|_{ L^1_T(\dot B^{\frac  dp}_{p,1})}+\|\q u\|^h_{ \widetilde{L}^\infty_T(\dot B^{-1+\frac  dp}_{p,1})}+\nu \|\q u\|^h_{L^1_T(\dot B^{1+\frac  dp}_{p,1})}
+\big\|(\q u_t+\nabla
 a)^h\big\|_{L^1_T(\dot{B}_{p,1}^{-1+\frac {d}{p}} )}\nonumber\\
&\quad\lesssim  \|\nu a_0\|^h_{\dot B^{\frac  dp}_{p,1}}+\|w_0\|^h_{\dot B^{-1+\frac  dp}_{p,1}}+(\v+\w)(\z+\v)+\nu^{-1}\y(\z+\v)\nonumber\\
&\quad\quad
+\nu^{-1}\x(\w+\v+\y)+\nu^{-2}\z\w+\nu^{-3}\x\y+\nu^{-2}\v^2\nonumber\\
&\quad\quad+\int_0^T\Big(\nu^{-1}\|a^h\|_{\dot B^{\frac  dp}_{p,1}}+\nu^{-1}\|\nu a^\ell\|_{\dot B^{1+\frac  d2}_{2,1}}\Big)\Big(\|\nu a^h\|_{\dot B^{\frac  dp}_{p,1}}+\| a^\ell\|_{\dot B^{-1+\frac  d2}_{2,1}}\Big)dt\nonumber\\
&\quad\quad
+\int_0^T\|(\nabla   u,\nabla V)\|_{\dot{B}_{p,1}^{\frac {d}{p}} }\Big(\|\q u^h\|_{\dot{B}_{p,1}^{-1+\frac {d}{p}}}+\|\q u^\ell\|_{\dot{B}_{2,1}^{-1+\frac {d}{2}} }+\|\nu a^h\|_{\dot B^{\frac  dp}_{p,1}}+\|\nu a^\ell\|_{\dot B^{\frac  d2}_{2,1}}\Big)dt.
\end{align}

\subsection{Low frequencies for the compressible part of \eqref{m}}

Now we  estimate the low frequencies for the compressible part  of \eqref{m} or (\ref{c}), which is rewritten as
\begin{eqnarray}\label{ccc}
\left\{\begin{aligned}
&a_t+(\p u+V)\cdot\nabla a+\diverg\q u  =-\diverg(a\,\q u),\\
&(\q u)_t-\nu\Delta\q u+\q( (u+V) \cdot \nabla \q u )+\nabla a=\q f,\\
&a|_{t=0}=a_0,\quad \q u|_{t=0}=\q v_0,
\end{aligned}\right.
\end{eqnarray}
with
\begin{align*}
f=&
a\left(\mathcal{V} _t+\p u_t+(\q u_t+\nabla a)\right)
+(1+a)(u+V)\cdot\nabla\p u \\
&+(1+a)(u+V)\cdot\nabla V+a(u+V)\cdot\nabla\q u +(k(a)-a)\nabla a.
\end{align*}
Applying $\ddj$ to both side of the first two equations of  \eqref{ccc}, we  get
\begin{align}
  &\partial_t\dot{\Delta}_ja +(\p u+\q u+V)\cdot\nabla \dot{\Delta}_ja +\diverg\q\dot{\Delta}_j u =\dot{\Delta}_jg,\label{p1a}\\
 &\partial_t\q\dot{\Delta}_j u+  \q( (u+V) \cdot \nabla\q\dot{\Delta}_j u )- \nu \Delta\q\dot{\Delta}_j u + \nabla \dot{\Delta}_ja = \dot{\Delta}_j\q f.\label{p1b}
 \end{align}
with
\begin{align*}
 \dot{\Delta}_jg :=&
\ddj(-a\diverg \q u)\!+\![(\p u\!+\!\q u\!+\!V)\cdot\nabla,\ddj] a,\\
\dot{\Delta}_jf:=&- [\ddj, u\!+\!V]\cdot\nabla \q u\!+\!\ddj\Big(a\left(\mathcal{V} _t\!+\!\p u_t\!+\!(\q u_t\!+\!\nabla a)\right)
\!+\!(1\!+\!a)\q u\cdot\nabla(\p u\!+\!V) \\
&\!+\!(1\!+\!a)(\p u\!+\!V)\cdot\nabla(\p u\!+\!V)\!+\!a(\p u\!+\!V)\cdot\nabla\q u \!+\!a\,\q u\cdot\nabla\q u\!+\!(k(a)-a)\nabla a\Big).
\end{align*}

We follow  from \cite{danchin2018} to bound each term $(\dot{\Delta}_ja,\q\dot{\Delta}_j u).$
More precisely, testing \eqref{p1a} and \eqref{p1b} by $\dot{\Delta}_ja$ and $\q\dot{\Delta}_j u,$ respectively, yields
\begin{align}\label{p5}
 \frac {1}{2} \frac {d}{dt} \int_{\R^d} (\dot{\Delta}_ja)^2\, dx + \int_{\R^d} \dot{\Delta}_ja\diverg\q\dot{\Delta}_j u\, dx =\frac {1}{2} \int_{\R^d} (\diverg \q u) \,  (\dot{\Delta}_ja)^2 \,dx+
  \int_{\R^d} \dot{\Delta}_jg \dot{\Delta}_ja\, dx
\end{align}
and
\begin{align}\label{p6}
\frac  12 \frac {d}{dt} \int_{\R^d} |\q\dot{\Delta}_j u|^2 \,dx + &\nu \int_{\R^d} |\nabla\q\dot{\Delta}_j u|^2 \,dx
- \int_{\R^d} \dot{\Delta}_ja \,\diverg\q\dot{\Delta}_j u \,dx\nonumber\\
&= \frac 12 \int_{\R^d} (\diverg \q u) |\q\dot{\Delta}_j u|^2\,dx +\int_{\R^d} \dot{\Delta}_j\q f\cdot\q\dot{\Delta}_j u\, dx.
\end{align}
Applying the gradient operator  $\nabla$ on both side of \eqref{p1a} gives
\begin{align}\label{p7}
 \nabla a_{j,t} +(\p u+\q u+V)\cdot\nabla\nabla \dot{\Delta}_ja + \nabla \diverg\q\dot{\Delta}_j u =
 \nabla \dot{\Delta}_jg-\nabla(\p u+\q u+V)\cdot\nabla \dot{\Delta}_ja.
\end{align}
Taking $L^2$ inner product  with $\nabla \dot{\Delta}_ja$ to the above equation  implies
\begin{align}\label{p8a}
&\frac 12 \frac {d}{dt}\int_{\R^d} |\nabla \dot{\Delta}_ja|^2 \, dx
+\int_{\R^d}\nabla\diverg\q\dot{\Delta}_j u\cdot\nabla \dot{\Delta}_ja\,dx\nonumber\\
&\quad=  \frac 12 \int_{\R^d} (\diverg \q u) |\nabla \dot{\Delta}_ja|^2\,dx+\int_{\R^d}(\nabla \dot{\Delta}_jg-\nabla(\p u+\q u+V)\cdot\nabla \dot{\Delta}_ja)\cdot\nabla \dot{\Delta}_ja\,dx.
\end{align}

The second term on the left hand side of (\ref{p8a}) is troublesome in our further estimation. We have to use some special technique to overcome the difficulty by
eliminating  this term, which involves the  highest order derivative. Actually,   it is convenient
to combine equation  (\ref{p8a}) with a relation involving $\int_{\R^d}\q\dot{\Delta}_j u\cdot \nabla \dot{\Delta}_ja\,dx.$ Now
testing  \eqref{p7}  by $\q\dot{\Delta}_j u$ and the  momentum equation (\ref{p1b}) by $\nabla \dot{\Delta}_ja,$ we get
\begin{align}\label{p8}
 &\frac {d}{dt}\int_{\R^d}\q\dot{\Delta}_j u \cdot\nabla \dot{\Delta}_ja\, dx+\int_{\R^d} (u+V)\cdot\nabla(\q\dot{\Delta}_j u\cdot\nabla \dot{\Delta}_ja)\,dx
 -\nu\int_{\R^d} \Delta\q\dot{\Delta}_j u \cdot\nabla a_{j}\,dx \nonumber\\
 &\quad\quad+ \int_{\R^d} |\nabla \dot{\Delta}_ja|^2 \,dx +\int_{\R^d}\nabla\diverg\q\dot{\Delta}_j u\cdot\q\dot{\Delta}_j u\,dx \nonumber\\
&\quad =
 \int_{\R^d} (\nabla \dot{\Delta}_jg-\nabla(\p u+\q u+V)\cdot\nabla \dot{\Delta}_ja) \cdot\q\dot{\Delta}_j u\,dx  +\int_{\R^d}  \dot{\Delta}_j\q f\cdot \nabla \dot{\Delta}_ja \,dx.
\end{align}

{Hence, multiplying  \eqref{p8a} by $\nu$ and adding the resultant equation to \eqref{p8}, we use the identity }  $\Delta\q\dot{\Delta}_j u \equiv \nabla \diverg\q\dot{\Delta}_j u$
to cancel the highest order terms to obtain  that
\begin{align}\label{p8-1}
&\frac 12 \frac {d}{dt}\int_{\R^d}\bigl(\nu|\nabla \dot{\Delta}_ja|^2+2\q\dot{\Delta}_j u \cdot\nabla \dot{\Delta}_ja\bigr)\,dx
+ \int_{\R^d}(|\nabla \dot{\Delta}_ja|^2-|\nabla\q\dot{\Delta}_j u|^2)\,dx\nonumber\\
&\quad=
\int_{\R^d} \biggl(\frac \nu2|\nabla \dot{\Delta}_ja|^2+\q\dot{\Delta}_j u\cdot\nabla \dot{\Delta}_ja\biggr)\diverg \q u\,dx
\!+\!\nu\int_{\R^d}\bigl(\nabla \dot{\Delta}_jg\!-\!  \nabla (u+V) \cdot\nabla \dot{\Delta}_ja)\cdot\nabla \dot{\Delta}_ja\,dx\nonumber\\
&\quad\quad+ \int_{\R^d} \bigl( \nabla \dot{\Delta}_jg -\nabla(u+V)\cdot\nabla \dot{\Delta}_ja\bigr)\cdot\q\dot{\Delta}_j u\,dx  +\int_{\R^d}  \dot{\Delta}_j\q f\cdot \nabla \dot{\Delta}_ja \,dx.
\end{align}
After multiplying (\ref{p8-1}) by $\nu$ and adding the resultant equation to  \eqref{p5} and \eqref{p6} respectively,  we get
\begin{align}\label{p9}
&\frac 12 \frac {d}{dt} L_j^2 +  \nu\int_{\R^d} \Bigl(|\nabla\q\dot{\Delta}_j u|^2 + |\nabla \dot{\Delta}_ja|^2\Bigr) dx \nonumber\\
&\quad=
 \int_{\R^d} \bigl(2\dot{\Delta}_jg\dot{\Delta}_ja+2\dot{\Delta}_j\q f\cdot\q\dot{\Delta}_j u+\nu^2\nabla \dot{\Delta}_jg\cdot\nabla \dot{\Delta}_ja+\nu\nabla \dot{\Delta}_jg\cdot\q\dot{\Delta}_j u+\nu \dot{\Delta}_j\q f\cdot\nabla \dot{\Delta}_ja\bigr)dx\nonumber\\
&\quad\quad +\frac 12\int_{\R^d}(2(\dot{\Delta}_ja)^2 + 2|\q\dot{\Delta}_j u|^2 + 2\nu\q\dot{\Delta}_j u \cdot\nabla \dot{\Delta}_ja + |\nu\nabla \dot{\Delta}_ja|^2)\diverg \q u\,dx\nonumber\\
&\quad\quad-\nu\int_{\R^d} \bigl(\nabla(u+V)\cdot\nabla \dot{\Delta}_ja\bigr)\cdot(\nu\nabla \dot{\Delta}_ja+\q\dot{\Delta}_j u)\,dx,
 \end{align}
where
\begin{equation}\label{p10}
 L_j := \left(\int_{\R^d}  \bigl(2(\dot{\Delta}_ja)^2 + 2|\q\dot{\Delta}_j u|^2 + 2\nu\q\dot{\Delta}_j u \cdot\nabla \dot{\Delta}_ja + |\nu\nabla \dot{\Delta}_ja|^2)\,dx\right)^{1/2}.
\end{equation}

It is readily seem  that
\begin{align}\label{p11}
C^{-1}\|(\q\dot{\Delta}_j u,\dot{\Delta}_ja,\nu\nabla \dot{\Delta}_ja)\|_{L^2}\le L_j\le C\|(\q\dot{\Delta}_j u,\dot{\Delta}_ja,\nu\nabla \dot{\Delta}_ja)\|_{L^2},
 \,\, j\in \Z,
\end{align}
and
\begin{align}\label{p11b}
 \nu\int ( |\nabla\q\dot{\Delta}_j u|^2 + |\nabla \dot{\Delta}_ja|^2)\, dx \geq C_2 \nu2^{2j}L_j^2,\,\, j\le j_0,
\end{align}
for  some constants $C, C_2>0$.
Therefore  it follows from \eqref{p9}, \eqref{p11} and \eqref{p11b}  that,  for $ j\le j_0$,
$$\frac 12\frac {d}{dt} L_j ^2+ \nu2^{2j} L_j^2\leq\Big(\frac 12\|\diverg \q u\|_{L^\infty}+C\|\nabla(u+V)\|_{L^\infty}\Big)L_j^2 + C\|[\dot{\Delta}_jg,\dot{\Delta}_j\q f,\nu\nabla \dot{\Delta}_jg]\|_{L^2}  L_j.$$
Hence, after integration in time, we have
\begin{align}
 L_j&(T)+ \nu2^{2j}\int_0^TL_j\,dt\nonumber\\
 &\leq L_j(0)+C\int_0^T\|\nabla(u+V)\|_{L^\infty}L_j\,dt
+ C\int_0^T\|[\dot{\Delta}_jg,\dot{\Delta}_j\q f,\nu\nabla \dot{\Delta}_jg]\|_{L^2}\,dt.\label{p13}
\end{align}
Summing up (\ref{p13}) with respect to  $j<j_0$, we have
\begin{align}\label{D00}
&\big\|(a^\ell,\nu \nabla a^\ell, \q u^\ell)\big\|_{L^\infty(0,T;\dot{B}_{2,1}^{-1+\frac {d}{2}} )}+\big\|(\nu a^\ell,\nu^2 \nabla a^\ell, \nu\q u^\ell)\big\|_{L^1(0,T;\dot{B}_{2,1}^{1+\frac {d}{2}} )}\nonumber\\
&\quad \le\big\|(a_0,\nu \nabla a_0, \q u_0)\big\|_{\dot{B}_{2,1}^{-1+\frac {d}{2}} }^\ell+ \int_0^T\|\nabla(u+V)\|_{L^\infty}\|(a^\ell,\nu\nabla a^\ell,\q u^\ell)\|_{\dot B^{-1+\frac  d2}_{2,1}}\,dt\nonumber\\
&\quad\quad+\int_0^T\sum_{j\in\Z}2^{(-1+\frac  d2)j}\| (\dot{\Delta}_jg,\nu\nabla \dot{\Delta}_jg)\|^\ell_{L^2}dt+\int_0^T\sum_{j\in\Z}2^{(-1+\frac  d2)j}\| \dot{\Delta}_j\q f\|^\ell_{L^2}dt.
\end{align}
In the following, we estimate last two terms on the right hand side of the previous  equality.

By Lemma \ref{xinqing}, we have
\begin{align}\label{D1}
&\int_0^T\sum_{j\in\Z}2^{(-1+\frac  d2)j}\| \ddj(-a\diverg \q u)\|^\ell_{L^2}dt\nonumber\\
&\quad\lesssim\int_0^T(\|\q u^\ell\|_{\dot{B}_{2,1}^{1+\frac {d}{2}}}+\|\q u^h\|_{\dot{B}_{p,1}^{1+\frac {d}{p}}})(\| a^\ell\|_{\dot{B}_{2,1}^{-1+\frac {d}{2}}}+\| a^h\|_{\dot{B}_{p,1}^{\frac {d}{p}}})dt\nonumber\\
&\quad\lesssim\int_0^T\nu^{-1}(\|\nu \q u^\ell\|_{\dot{B}_{2,1}^{1+\frac {d}{2}}}+\|\nu \q u^h\|_{\dot{B}_{p,1}^{1+\frac {d}{p}}})\| a^\ell\|_{\dot{B}_{2,1}^{-1+\frac {d}{2}}}dt\nonumber\\
&\quad\quad+\int_0^T\nu^{-2}(\|\nu \q u^\ell\|_{\dot{B}_{2,1}^{1+\frac {d}{2}}}+\|\nu \q u^h\|_{\dot{B}_{p,1}^{1+\frac {d}{p}}})\|\nu a^h\|_{\dot{B}_{p,1}^{\frac {d}{p}}}dt.
\end{align}
By Lemma \ref{xinqing}, we have
\begin{align}\label{D2}
&\int_0^T\sum_{j\in\Z}2^{(-1+\frac  d2)j}\|[(\p u+\q u+V)\cdot\nabla,\ddj]  a\|^\ell_{L^2}dt\nonumber\\
&\quad\lesssim(\|(\nabla \p u,\nabla V)\|_{\dot{B}_{p,1}^{\frac {d}{p}}}+\|\nabla \q u^h\|_{\dot{B}_{p,1}^{\frac {d}{p}}}+\|\nabla \q u^\ell\|_{\dot{B}_{2,1}^{\frac {d}{2}}})(\| a^\ell\|_{\dot{B}_{2,1}^{-1+\frac {d}{2}}}+\| a^h\|_{\dot{B}_{p,1}^{\frac {d}{p}}})\nonumber\\
&\quad\lesssim\int_0^T\Big(\nu^{-1}(\|\nu \q u^\ell\|_{\dot{B}_{2,1}^{1+\frac {d}{2}}}+\|\nu \q u^h\|_{\dot{B}_{p,1}^{1+\frac {d}{p}}})+\|(\nabla \p u,\nabla V)\|_{\dot{B}_{p,1}^{\frac {d}{p}}}\Big)\| a^\ell\|_{\dot{B}_{2,1}^{-1+\frac {d}{2}}}dt\nonumber\\
&\quad\quad+\int_0^T\Big(\nu^{-2}(\|\nu \q u^\ell\|_{\dot{B}_{2,1}^{1+\frac {d}{2}}}+\|\nu \q u^h\|_{\dot{B}_{p,1}^{1+\frac {d}{p}}})+\nu^{-1}\|(\nabla \p u,\nabla V)\|_{\dot{B}_{p,1}^{\frac {d}{p}}}\Big)\|\nu a^h\|_{\dot{B}_{p,1}^{\frac {d}{p}}}dt.
\end{align}
From the above two estimates, we  find that
\begin{align}\label{D3}
&\int_0^T\sum_{j\in\Z}2^{(-1+\frac  d2)j}\| \dot{\Delta}_jg\|^\ell_{L^2}dt
\nonumber\\
&\quad\lesssim\int_0^T\Big(\nu^{-1}(\|\nu \q u^\ell\|_{\dot{B}_{2,1}^{1+\frac {d}{2}}}+\|\nu \q u^h\|_{\dot{B}_{p,1}^{1+\frac {d}{p}}})+\|(\nabla \p u,\nabla V)\|_{\dot{B}_{p,1}^{\frac {d}{p}}}\Big)\| a^\ell\|_{\dot{B}_{2,1}^{-1+\frac {d}{2}}}dt\nonumber\\
&\quad\quad+\int_0^T\Big(\nu^{-2}(\|\nu \q u^\ell\|_{\dot{B}_{2,1}^{1+\frac {d}{2}}}+\|\nu \q u^h\|_{\dot{B}_{p,1}^{1+\frac {d}{p}}})+\nu^{-1}\|(\nabla \p u,\nabla V)\|_{\dot{B}_{p,1}^{\frac {d}{p}}}\Big)\|\nu a^h\|_{\dot{B}_{p,1}^{\frac {d}{p}}}dt.
\end{align}
Similarly, from Lemma \ref{xinqing}, we  get
\begin{align}\label{D5}
&\nu\int_0^T\sum_{j\in\Z}2^{(-1+\frac  d2)j}\| \nabla\ddj(-a\diverg \q u)\|^\ell_{L^2}dt\nonumber\\
&\quad\lesssim\int_0^T(\|\q u^\ell\|_{\dot{B}_{2,1}^{1+\frac {d}{2}}}+\|\q u^h\|_{\dot{B}_{p,1}^{1+\frac {d}{p}}})(\|\nu a^\ell\|_{\dot{B}_{2,1}^{\frac {d}{2}}}+\|\nu a^h\|_{\dot{B}_{p,1}^{\frac {d}{p}}})dt\nonumber\\
&\quad\lesssim\int_0^T\nu^{-1}(\|\nu\nabla \q u^h\|_{\dot{B}_{p,1}^{\frac {d}{p}}}+\|\nu\nabla \q u^\ell\|_{\dot{B}_{2,1}^{\frac {d}{2}}})(\|\nu a^\ell\|_{\dot{B}_{2,1}^{\frac {d}{2}}}+\|\nu a^h\|_{\dot{B}_{p,1}^{\frac {d}{p}}})dt,
\end{align}
and
\begin{align}\label{D6}
&\int_0^T\sum_{j\in\Z}2^{(-1+\frac  d2)j}\|\nabla([(\p u+\q u+V)\cdot\nabla,\ddj] \nu a)\|^\ell_{L^2}dt\nonumber\\
&\quad\lesssim\int_0^T\|(\nabla \p u,\nabla V)\|_{\dot{B}_{p,1}^{\frac {d}{p}}}(\|\nu a^\ell\|_{\dot{B}_{2,1}^{\frac {d}{2}}}+\|\nu a^h\|_{\dot{B}_{p,1}^{\frac {d}{p}}})dt\nonumber\\
&\quad\quad+\int_0^T\nu^{-1}(\|\nu\nabla \q u^h\|_{\dot{B}_{p,1}^{\frac {d}{p}}}+\|\nu\nabla \q u^\ell\|_{\dot{B}_{2,1}^{\frac {d}{2}}})(\|\nu a^\ell\|_{\dot{B}_{2,1}^{\frac {d}{2}}}+\|\nu a^h\|_{\dot{B}_{p,1}^{\frac {d}{p}}})dt.
\end{align}

 Hence
\begin{align}\label{D8}
&\nu\int_0^T\sum_{j\in\Z}2^{(-1+\frac  d2)j}\|\nabla \dot{\Delta}_jg\|^\ell_{L^2}dt
\nonumber\\
&\quad\lesssim\int_0^T\|(\nabla \p u,\nabla V)\|_{\dot{B}_{p,1}^{\frac {d}{p}}}(\|\nu a^\ell\|_{\dot{B}_{2,1}^{\frac {d}{2}}}+\|\nu a^h\|_{\dot{B}_{p,1}^{\frac {d}{p}}})dt\nonumber\\
&\quad\quad+\int_0^T\nu^{-1}(\|\nu\nabla \q u^h\|_{\dot{B}_{p,1}^{\frac {d}{p}}}+\|\nu\nabla \q u^\ell\|_{\dot{B}_{2,1}^{\frac {d}{2}}})(\|\nu a^\ell\|_{\dot{B}_{2,1}^{\frac {d}{2}}}+\|\nu a^h\|_{\dot{B}_{p,1}^{\frac {d}{p}}})dt.
\end{align}

Now  we claim that, for any $p$  given  in Theorem \ref{dingli}, there holds
\begin{align}\label{pengyou}
\|\q(bc)\|_{\dot{B}_{2,1}^{-1+\frac {d}{2}}}^\ell\lesssim(\|b\|_{\dot{B}_{p,1}^{-1+\frac {d}{p}}}+\|b^\ell\|_{\dot{B}_{2,1}^{-1+\frac {d}{2}}})\|c\|_{\dot{B}_{p,1}^{\frac {d}{p}}}.
\end{align}
Indeed,
 denoting $\q^\ell:=\dot S_{j_0+1}\q,$ we  get
\begin{align}\label{D9}
\q^\ell(bc)=\q^\ell\bigl(\dot{T}_{b} c+ \dot{R}(b,c)\bigr)+\dot{T}_{c}\q^\ell b
+[\q^\ell,T_{c}]b.
\end{align}
By Lemma \ref{fenjie}, we obtain
\begin{align}\label{D11-1}
\|\dot{T}_{c}\q^\ell b\|_{\dot{B}_{2,1}^{-1+\frac {d}{2}}}\lesssim\|c\|_{L^\infty}\|b^\ell\|_{\dot{B}_{2,1}^{-1+\frac {d}{2}}}
\lesssim&\|b^\ell\|_{\dot{B}_{2,1}^{-1+\frac {d}{2}}}\|c\|_{\dot{B}_{p,1}^{\frac {d}{p}}},\nonumber\\
\|\q^\ell\bigl(T_{b} c+ \dot{R}(b,c)\bigr)\|_{\dot{B}_{2,1}^{-1+\frac {d}{2}}}\lesssim\|b\|_{\dot{B}_{p^*,1}^{-1+\frac {d}{p^*}}}\|c\|_{\dot{B}_{p,1}^{\frac {d}{p}}}
\lesssim&\|b\|_{\dot{B}_{p,1}^{-1+\frac {d}{p}}}\|c\|_{\dot{B}_{p,1}^{\frac {d}{p}}}.
\end{align}
By Lemma \ref{dahai},   we have
\begin{align}\label{D11}
\|[\q^\ell,\dot{T}_{c}]b\|_{\dot{B}_{2,1}^{-1+\frac {d}{2}}}\lesssim\|\nabla c\|_{\dot{B}_{p^*,1}^{-1+\frac {d}{p^*}}}\|b\|_{\dot{B}_{p,1}^{-1+\frac {d}{p}}}
\lesssim&\|b\|_{\dot{B}_{p,1}^{-1+\frac {d}{p}}}\|c\|_{\dot{B}_{p,1}^{\frac {d}{p}}}, \,\,\,\,\,\frac  1p+\frac  1{p*}=\frac  12.
\end{align}
Thus, the combination of (\ref{D9})--(\ref{D11}) shows the validity of  \eqref{pengyou}.

Moreover, if $\diverg b=0,$ then we have  $\q^\ell b=0$ and thus
\begin{align}\label{pengyou2}
\|\q(bc)\|_{\dot{B}_{2,1}^{-1+\frac {d}{2}}}^\ell\lesssim\|b\|_{\dot{B}_{p,1}^{-1+\frac {d}{p}}}\|c\|_{\dot{B}_{p,1}^{\frac {d}{p}}}.
\end{align}
Since $ \diverg ({V} _t+\p u_t)=0$,
taking  $b={V} _t+\p u_t$ and $c= a$ in \eqref{pengyou2}, we have
\begin{align}\label{D12}
&\|\q\left( a({V} _t+\p u_t)\right)\|_{L^1(0,T;\dot{B}_{2,1}^{-1+\frac {d}{2}} )}^\ell\lesssim \|a\|_{L^\infty(0,T;\dot{B}_{p,1}^{\frac {d}{p}} )}\|(\mathcal{V} _t,\p u_t)\|_{L^1(0,T;\dot{B}_{p,1}^{-1+\frac {d}{p}} )}\lesssim\nu^{-1}\x(\v+\w).
\end{align}
Similarly, using $\|a\|_{L^\infty(0,T;\dot B^{\frac  d2}_{2,1})}\ll1,$ we have
\begin{equation}\label{D13}
\|\q\left( (1+a)(\p u+{V} )\cdot\nabla (\p u+{V})\right)\|_{L^1(0,T;\dot{B}_{2,1}^{-1+\frac {d}{2}} )}^\ell
\lesssim(\v+\w)(\z+\v),
 \end{equation}
\begin{equation}\label{D15}
\|\q\left( a(\p u+{V} )\cdot\nabla \q u\right)\|_{L^1(0,T;\dot{B}_{2,1}^{-1+\frac {d}{2}} )}^\ell
\lesssim\nu^{-2}\x\y(\z+\v).
 \end{equation}
 Taking $b=\q u_t+\nabla a$ and $c= a$ in  \eqref{pengyou} gives
\begin{align}\label{D16}
&\|\q\left( a(\q u_t+\nabla a)\right)\|_{L^1(0,T;\dot{B}_{2,1}^{-1+\frac {d}{2}} )}^\ell\nonumber\\
&\quad\lesssim\|a\|_{L^\infty(0,T;\dot{B}_{p,1}^{\frac {d}{p}} )}\Big(\|(\q u_t+\nabla a)^h\|_{L^1(0,T;\dot{B}_{p,1}^{-1+\frac {d}{p}} )}
+\|(\q u_t+\nabla a)^\ell\|_{L^1(0,T;\dot{B}_{2,1}^{-1+\frac {d}{2}} )}
\Big)\nonumber\\
&\quad\lesssim\nu^{-1}\x\y.
 \end{align}
Similarly, we obtain that
\begin{align}\label{D18}
&\|\q\left( (1+a)\q u\cdot\nabla (\p u+{V})\right)\|_{L^1(0,T;\dot{B}_{2,1}^{-1+\frac {d}{2}} )}^\ell
\nonumber\\
&\quad\lesssim\int_0^T\|(\nabla  \p u,\nabla V)\|_{\dot{B}_{p,1}^{\frac {d}{p}} }
\Big(\|\q u^h\|_{\dot{B}_{p,1}^{-1+\frac {d}{p}})}+\|\q u^\ell\|_{\dot{B}_{2,1}^{-1+\frac {d}{2}} }\Big)dt,
 \end{align}
 and
\begin{align}\label{D19}
&\|\q\left( a\q u\cdot\nabla \q u\right)\|_{L^1(0,T;\dot{B}_{2,1}^{-1+\frac {d}{2}} )}^\ell
\lesssim\nu^{-2}\x^2\y.
 \end{align}
By Proposition \ref{fuhe2}, we have
\begin{align}\label{D20}
&\|\q\left((k(a)-a)\nabla a\right)\|_{L^1(0,T;\dot{B}_{2,1}^{-1+\frac {d}{2}} )}^\ell\nonumber\\
&\quad\lesssim\int_0^T
(\|\nabla a\|_{\dot{B}_{p,1}^{-1+\frac {d}{p}}}+\|\nabla a^\ell\|_{\dot{B}_{2,1}^{-1+\frac {d}{2}}})\|a\|_{\dot{B}_{p,1}^{\frac {d}{p}}}dt
\nonumber\\
&\quad\lesssim\|a\|_{L^2(0,T;\dot{B}_{p,1}^{\frac {d}{p}} )}^2+\|a\|_{L^2(0,T;\dot{B}_{p,1}^{\frac {d}{p}} )}\|a^\ell\|_{L^2(0,T;\dot{B}_{2,1}^{\frac {d}{2}} )}\nonumber\\
&\quad\lesssim\|a^h\|_{L^\infty(0,T;\dot{B}_{p,1}^{\frac {d}{p}} )}\|a^h\|_{L^1(0,T;\dot{B}_{p,1}^{\frac {d}{p}} )}+\|a^\ell\|_{L^\infty(0,T;\dot{B}_{2,1}^{-1+\frac {d}{2}} )}\|a^\ell\|_{L^1(0,T;\dot{B}_{2,1}^{1+\frac {d}{2}} )}\nonumber\\
&\quad\lesssim\int_0^T\Big(\nu^{-1}\|a^h\|_{\dot{B}_{p,1}^{\frac {d}{p}} }+\nu^{-1}\|\nu a^\ell\|_{\dot{B}_{2,1}^{1+\frac {d}{2}} }\Big)\Big(\|\nu a^h\|_{\dot{B}_{p,1}^{\frac {d}{p}} }+\| a^\ell\|_{\dot{B}_{2,1}^{-1+\frac {d}{2}} }\Big)dt.
 \end{align}
It follows from  Lemma \ref{xinqing} that
\begin{align}\label{D22}
&\int_0^T\sum_{j\in\Z}2^{(-1+\frac  d2)j}\|[\ddj, u+V]\cdot\nabla \q u\|^\ell_{L^2}dt\nonumber\\
&\quad\lesssim\int_0^T(\|(\nabla \p u,\nabla \q u^h,\nabla V)\|_{\dot{B}_{p,1}^{\frac {d}{p}}}+\|\nabla \q u^\ell\|_{\dot{B}_{2,1}^{\frac {d}{2}}})(\|\q u^\ell\|_{\dot{B}_{2,1}^{-1+\frac {d}{2}}}+\|\q u^h\|_{\dot{B}_{p,1}^{-1+\frac {d}{p}}})dt.
\end{align}

Now we assume $
\nu^{-1}\ll1$.
Inserting  \eqref{D3}, \eqref{D8}, \eqref{D12}--\eqref{D22} into  \eqref{D00} gives
\begin{align}\label{D23}
&\big\|(a^\ell,\nu \nabla a^\ell, \q u^\ell)\big\|_{L^\infty(0,T;\dot{B}_{2,1}^{-1+\frac {d}{2}} )}+\big\|(\nu  a^\ell,\nu^2 \nabla a^\ell, \nu\q u^\ell)\big\|_{L^1(0,T;\dot{B}_{2,1}^{1+\frac {d}{2}} )}\nonumber\\
&\quad \lesssim\big\|(a_0,\nu \nabla a_0, \q u_0)\big\|_{\dot{B}_{2,1}^{-1+\frac {d}{2}} }^\ell+\int_0^T\|(\nabla \p u,\nabla V)\|_{\dot{B}_{p,1}^{\frac {d}{p}}}(\| a^\ell\|_{\dot{B}_{2,1}^{-1+\frac {d}{2}}}+\|\nu a^\ell\|_{\dot{B}_{2,1}^{\frac {d}{2}}}+\|\nu a^h\|_{\dot{B}_{p,1}^{\frac {d}{p}}})dt\nonumber\\
&\quad\quad+\int_0^T\Big(\nu^{-1}\|a^h\|_{\dot{B}_{p,1}^{\frac {d}{p}} }+\nu^{-1}\|\nu a^\ell\|_{\dot{B}_{2,1}^{1+\frac {d}{2}} }\Big)(\| a^\ell\|_{\dot{B}_{2,1}^{-1+\frac {d}{2}}}+\|\nu a^h\|_{\dot{B}_{p,1}^{\frac {d}{p}}})dt
\nonumber\\
&\quad\quad+\int_0^T\Big(\|\q u^\ell\|_{\dot{B}_{2,1}^{1+\frac {d}{2}}}+\|\q u^h\|_{\dot{B}_{p,1}^{1+\frac {d}{p}}}\Big)(\| a^\ell\|_{\dot{B}_{2,1}^{-1+\frac {d}{2}}}+\|\nu a^\ell\|_{\dot{B}_{2,1}^{\frac {d}{2}}}+\|\nu a^h\|_{\dot{B}_{p,1}^{\frac {d}{p}}})dt\nonumber\\
&\quad\quad+\int_0^T\Big(\|(\nabla \p u,\nabla \q u^h,\nabla V)\|_{\dot{B}_{p,1}^{\frac {d}{p}}}+\|\nabla \q u^\ell\|_{\dot{B}_{2,1}^{\frac {d}{2}}}\Big)
\Big(\|\q u^h\|_{\dot{B}_{p,1}^{-1+\frac {d}{p}})}+\|\q u^\ell\|_{\dot{B}_{2,1}^{-1+\frac {d}{2}} }\Big)dt\nonumber\\
&\quad\quad+
\nu^{-1}\x(\v+\w)+\nu^{-1}\x\y+(\v+\w)(\z+\v)+\nu^{-2}\x\y(\z+\v)+\nu^{-2}\x^2\y.
\end{align}
Form the second equation in \eqref{ccc}, we deduce that
\begin{align}\label{D2326}
&\big\|(\q u_t+\nabla a)^\ell\big\|_{L^1(0,T;\dot{B}_{2,1}^{-1+\frac {d}{2}} )}\nonumber\\
&\quad\lesssim\|\nu\q u\|^\ell_{L^1(0,T;\dot{B}_{2,1}^{-1+\frac {d}{2}} )}+\|\q( (u+V) \cdot \nabla \q u )\|^\ell_{L^1(0,T;\dot{B}_{2,1}^{-1+\frac {d}{2}} )}+\|\q f\|^\ell_{L^1(0,T;\dot{B}_{2,1}^{-1+\frac {d}{2}} )}.
\end{align}
Following the  derivation of \eqref{D15} and \eqref{D19}, we  get
\begin{align}\label{D2323}
&\|\q( (u+V) \cdot \nabla \q u )\|^\ell_{L^1(0,T;\dot{B}_{2,1}^{-1+\frac {d}{2}} )}\nonumber\\
&\quad\lesssim\Big(\|\q u^h\|_{L^\infty(0,T;\dot{B}_{p,1}^{-1+\frac {d}{p}} )}+\|\q u^\ell\|_{L^\infty(0,T;\dot{B}_{2,1}^{-1+\frac {d}{2}} )}+\|( \p u, V)\|_{L^\infty(0,T;\dot{B}_{p,1}^{-1+\frac {d}{p}} )}\Big)\|\nabla \q u\|_{L^1(0,T;\dot{B}_{p,1}^{\frac {d}{p}} )}\nonumber\\
&\quad\lesssim\nu^{-1}\y(\x+\z+\v).
 \end{align}
 Combining  the estimates \eqref{D12}--\eqref{D22}, we find that
\begin{align}\label{D22+963}
&\|\q f\|^\ell_{L^1(0,T;\dot{B}_{2,1}^{-1+\frac {d}{2}} )}
\lesssim\int_0^T\Big(\nu^{-1}\|a^h\|_{\dot{B}_{p,1}^{\frac {d}{p}} }+\nu^{-1}\|\nu a^\ell\|_{\dot{B}_{2,1}^{1+\frac {d}{2}} }\Big)(\| a^\ell\|_{\dot{B}_{2,1}^{-1+\frac {d}{2}}}+\|\nu a^h\|_{\dot{B}_{p,1}^{\frac {d}{p}}})dt\nonumber\\
&\quad\quad+\int_0^T\Big(\|(\nabla \p u,\nabla \q u^h,\nabla V)\|_{\dot{B}_{p,1}^{\frac {d}{p}}}+\|\nabla \q u^\ell\|_{\dot{B}_{2,1}^{\frac {d}{2}}}\Big)
\Big(\|\q u^h\|_{\dot{B}_{p,1}^{-1+\frac {d}{p}})}+\|\q u^\ell\|_{\dot{B}_{2,1}^{-1+\frac {d}{2}} }\Big)dt\nonumber\\
&\quad\quad+
\nu^{-1}\x(\v+\w)+\nu^{-1}\x\y+(\v+\w)(\z+\v)+\nu^{-2}\x\y(\z+\v)+\nu^{-2}\x^2\y.
\end{align}

By  \eqref{D23}--\eqref{D22+963}, we have
\begin{align}\label{DD23}
&\big\|(a^\ell,\nu \nabla a^\ell, \q u^\ell)\big\|_{L^\infty(0,T;\dot{B}_{2,1}^{-1+\frac {d}{2}} )}+\big\|(\nu  a^\ell,\nu^2 \nabla a^\ell, \nu\q u^\ell)\big\|_{L^1(0,T;\dot{B}_{2,1}^{1+\frac {d}{2}} )}\nonumber\\
&\quad\quad+\big\|(\q u_t+\nabla a)^\ell\big\|_{L^1(0,T;\dot{B}_{2,1}^{-1+\frac {d}{2}} )}\nonumber\\
&\quad \lesssim\big\|(a_0,\nu \nabla a_0, \q u_0)\big\|_{\dot{B}_{2,1}^{-1+\frac {d}{2}} }^\ell+
\nu^{-1}\x(\v+\w)+\nu^{-1}\y(\x+\z+\v)
\nonumber\\
&\quad\quad+(\v+\w)(\z+\v)+\nu^{-2}\x\y(\z+\v)+\nu^{-2}\x^2\y\nonumber\\
&\quad\quad+\int_0^T\Big(\nu^{-1}\|a^h\|_{\dot{B}_{p,1}^{\frac {d}{p}} }+\nu^{-1}\|\nu a^\ell\|_{\dot{B}_{2,1}^{1+\frac {d}{2}} }\Big)(\| a^\ell\|_{\dot{B}_{2,1}^{-1+\frac {d}{2}}}+\|\nu a^h\|_{\dot{B}_{p,1}^{\frac {d}{p}}})dt\nonumber\\
&\quad\quad
+\int_0^T\|(\nabla \p u,\nabla V)\|_{\dot{B}_{p,1}^{\frac {d}{p}}}(\| a^\ell\|_{\dot{B}_{2,1}^{-1+\frac {d}{2}}}+\|\nu a^\ell\|_{\dot{B}_{2,1}^{\frac {d}{2}}}+\|\nu a^h\|_{\dot{B}_{p,1}^{\frac {d}{p}}})dt\nonumber\\
&\quad\quad+\int_0^T\Big(\|\q u^\ell\|_{\dot{B}_{2,1}^{1+\frac {d}{2}}}+\|\q u^h\|_{\dot{B}_{p,1}^{1+\frac {d}{p}}}\Big)(\| a^\ell\|_{\dot{B}_{2,1}^{-1+\frac {d}{2}}}+\|\nu a^\ell\|_{\dot{B}_{2,1}^{\frac {d}{2}}}+\|\nu a^h\|_{\dot{B}_{p,1}^{\frac {d}{p}}})dt\nonumber\\
&\quad\quad+\int_0^T\Big(\|(\nabla \p u,\nabla \q u^h,\nabla V)\|_{\dot{B}_{p,1}^{\frac {d}{p}}}+\|\nabla \q u^\ell\|_{\dot{B}_{2,1}^{\frac {d}{2}}}\Big)
\Big(\|\q u^h\|_{\dot{B}_{p,1}^{-1+\frac {d}{p}})}+\|\q u^\ell\|_{\dot{B}_{2,1}^{-1+\frac {d}{2}} }\Big)dt.
\end{align}

By the  estimates \eqref{W123} and  \eqref{DD23},
we  get
\begin{align}\label{D25}
&\big\|(a^\ell,\nu \nabla a^\ell, \q u^\ell)\big\|_{\widetilde{L}^\infty_T(\dot{B}_{2,1}^{-1+\frac {d}{2}} )}+\big\|(\nu  a^\ell,\nu^2 \nabla a^\ell, \nu\q u^\ell)\big\|_{L^1_T(\dot{B}_{2,1}^{1+\frac {d}{2}} )}\nonumber\\
&\quad\quad+\|\nu a\|^h_{\widetilde{L}^\infty_T(\dot B^{\frac  dp}_{p,1})}+\|a^h\|_{ L^1_T(\dot B^{\frac  dp}_{p,1})}+\|\q u\|^h_{ \widetilde{L}^\infty_T(\dot B^{-1+\frac  dp}_{p,1})}+\nu \|\q u\|^h_{L^1_T(\dot B^{1+\frac  dp}_{p,1})}\nonumber\\
&\quad\quad
+\big\|(\q u_t+\nabla a)^\ell\big\|_{L^1(0,T;\dot{B}_{2,1}^{-1+\frac {d}{2}} )}+\big\|(\q u_t+\nabla a)^h\big\|_{L^1(0,T;\dot{B}_{p,1}^{-1+\frac {d}{p}} )}
\nonumber\\
&\quad \lesssim\big\|(a_0,\nu \nabla a_0, \q u_0)\big\|_{\dot{B}_{2,1}^{-1+\frac {d}{2}} }^\ell+\|\nu a_0\|^h_{\dot B^{\frac  dp}_{p,1}}+\|\nu^{-1} a_0\|^h_{\dot B^{\frac  dp}_{p,1}}+\|\q u_0\|^h_{\dot B^{-1+\frac  dp}_{p,1}}\nonumber\\
&\quad\quad
+(\v+\w)(\z+\v)+\nu^{-2}\x\y(\z+\v)+\nu^{-2}\x^2\y+\nu^{-1}\y(\z+\v)\nonumber\\
&\quad\quad
+\nu^{-1}\x(\w+\v+\y)+\nu^{-2}\z\w+\nu^{-3}\x\y+\nu^{-2}\v^2\nonumber\\
&\quad\quad
+\int_0^T\Big(\nu^{-1}\|a^h\|_{\dot{B}_{p,1}^{\frac {d}{p}} }+\nu^{-1}\|\nu a^\ell\|_{\dot{B}_{2,1}^{1+\frac {d}{2}} }\Big)(\| a^\ell\|_{\dot{B}_{2,1}^{-1+\frac {d}{2}}}+\|\nu a^h\|_{\dot{B}_{p,1}^{\frac {d}{p}}})dt\nonumber\\
&\quad\quad
+\int_0^T\|(\nabla \p u,\nabla V)\|_{\dot{B}_{p,1}^{\frac {d}{p}}}(\| a^\ell\|_{\dot{B}_{2,1}^{-1+\frac {d}{2}}}+\|\nu a^\ell\|_{\dot{B}_{2,1}^{\frac {d}{2}}}+\|\nu a^h\|_{\dot{B}_{p,1}^{\frac {d}{p}}})dt\nonumber\\
&\quad\quad+\int_0^T\Big(\|\q u^\ell\|_{\dot{B}_{2,1}^{1+\frac {d}{2}}}+\|\q u^h\|_{\dot{B}_{p,1}^{1+\frac {d}{p}}}\Big)(\| a^\ell\|_{\dot{B}_{2,1}^{-1+\frac {d}{2}}}+\|\nu a^\ell\|_{\dot{B}_{2,1}^{\frac {d}{2}}}+\|\nu a^h\|_{\dot{B}_{p,1}^{\frac {d}{p}}})dt\nonumber\\
&\quad\quad+\int_0^T\Big(\|(\nabla \p u,\nabla \q u^h,\nabla V)\|_{\dot{B}_{p,1}^{\frac {d}{p}}}+\|\nabla \q u^\ell\|_{\dot{B}_{2,1}^{\frac {d}{2}}}\Big)
\Big(\|\q u^h\|_{\dot{B}_{p,1}^{-1+\frac {d}{p}})}+\|\q u^\ell\|_{\dot{B}_{2,1}^{-1+\frac {d}{2}} }\Big)dt.
\end{align}
Hence,  from the Gronwall inequality, we have
\begin{align}\label{D26}
\x+\y\le&\exp\Bigg\{C\int_0^T\Big(\|(\nabla \p u,\nabla \q u^h,\nabla V)\|_{\dot{B}_{p,1}^{\frac {d}{p}}}+\|\nabla \q u^\ell\|_{\dot{B}_{2,1}^{\frac {d}{2}}}+\nu^{-1}\|a^h\|_{\dot{B}_{p,1}^{\frac {d}{p}} }\nonumber\\
&+\nu^{-1}\|\nu a^\ell\|_{\dot{B}_{2,1}^{1+\frac {d}{2}} }\Big)dt\Bigg\}
\Bigg\{\x(0)+(\v+\w)(\z+\v)+\nu^{-2}\x\y(\z+\v)+\nu^{-2}\x^2\y\nonumber\\
&+\nu^{-1}\y(\z+\v)
+\nu^{-1}\x(\w+\v+\y)+\nu^{-2}\z\w+\nu^{-3}\x\y+\nu^{-2}\v^2\Bigg\}.
\end{align}

\subsection{ Existence of the global-in-time solution}
We deduce from \eqref{D26} that
\begin{align}\label{T3}
\x+\y\le &\exp\Big(C(\w+\v+\nu^{-1}\y)\Big)
\Bigg\{\x(0)+(\v+\w)(\z+\v)+\nu^{-2}\x\y(\z+\v)\nonumber\\
&\hspace{-0.3cm} +\nu^{-2}\x^2\y+\nu^{-1}\y(\z+\v)
+\nu^{-1}\x(\w+\v+\y)+\nu^{-2}\z\w+\nu^{-3}\x\y+\nu^{-2}\v^2\Bigg\}.
\end{align}
%
This together with  the assumption $\nu^{-1}\ll 1$ and \eqref{T1} yields
\begin{align}\label{T6}
\x+\y\le &\exp\Big(C(M+\delta+\nu^{-1}{\eta})\Big)
\Bigg\{\x(0)+(M+\delta)^2+\nu^{-1}{\eta}(\delta+M)\x+\nu^{-1}{\eta}^2\x\nonumber\\
&\hspace{5cm}
+\nu^{-1}(M+\delta+{\eta})\x+\nu^{-1}(M+\delta)\y\Bigg\}.
\end{align}
Recalling that
\begin{align}\label{T8}
\z+\w
\le&\exp\Big(C\int_0^T\Big(\|(\nabla V,\nabla \q u)\|_{\dot{B}_{p,1}^{\frac {d}{p}} }+\|u\|_{\dot{B}_{p,1}^{\frac {d}{p}} }^2+\|V\|_{\dot{B}_{p,1}^{\frac {d}{p}} }^2\Big)dt\Big)\Bigg\{\nu^{-\frac 12}\x^{\frac 12}\y^{\frac 12}\v\nonumber\\
&\quad\quad\quad\quad+\nu^{-1} \x(\z+\x+\v)\w+\nu^{-1}(\y+\w+\v)\x+\nu^{-1} \x\v^2\Bigg\},
\end{align}
and employing  an embedding inequality and the interpolation inequality, we have
\begin{align*}
\|u\|_{\dot{B}_{p,1}^{\frac {d}{p}} }^2+\|V\|_{\dot{B}_{p,1}^{\frac {d}{p}} }^2
\lesssim&\|\p u\|_{\dot{B}_{p,1}^{\frac {d}{p}} }^2+\|\q u^h\|_{\dot{B}_{p,1}^{\frac {d}{p}} }^2+\|\q u^\ell\|_{\dot{B}_{2,1}^{\frac {d}{2}} }^2+\|V\|_{\dot{B}_{p,1}^{\frac {d}{p}} }^2\nonumber\\
\lesssim&\|\p u\|_{\dot{B}_{p,1}^{-1+\frac {d}{p}} }\|\p u\|_{\dot{B}_{p,1}^{1+\frac {d}{p}} }+\nu^{-1}\|\q u^h\|_{\dot{B}_{p,1}^{-1+\frac {d}{p}} }\|\nu\q u^h\|_{\dot{B}_{p,1}^{1+\frac {d}{p}} }\nonumber\\
&\quad+\nu^{-1}\|\q u^\ell\|_{\dot{B}_{2,1}^{-1+\frac {d}{2}} }\|\nu\q u^\ell\|_{\dot{B}_{2,1}^{1+\frac {d}{2}} }+\|V\|_{\dot{B}_{p,1}^{-1+\frac {d}{p}} }\|V\|_{\dot{B}_{p,1}^{1+\frac {d}{p}} },
\end{align*}
and hence, by \eqref{T1},
\begin{align}\label{TT11}
&\int_0^T\Big(\|(\nabla V,\nabla \q u)\|_{\dot{B}_{p,1}^{\frac {d}{p}} }+\|u\|_{\dot{B}_{p,1}^{\frac {d}{p}} }^2+\|V\|_{\dot{B}_{p,1}^{\frac {d}{p}} }^2\Big)dt\nonumber\\
&\lesssim\int_0^T
\Big(\|(\nabla V,\nabla \q u^h)\|_{\dot{B}_{p,1}^{\frac {d}{p}}}+\|\nabla \q u^\ell\|_{\dot{B}_{2,1}^{\frac {d}{2}}}+\nu^{-1}\|\q u^h\|_{\dot{B}_{p,1}^{-1+\frac {d}{p}} }\|\nu\q u^h\|_{\dot{B}_{p,1}^{1+\frac {d}{p}} }\nonumber\\
&\quad\quad+\|\p u\|_{\dot{B}_{p,1}^{-1+\frac {d}{p}} }\|\p u\|_{\dot{B}_{p,1}^{1+\frac {d}{p}} }+\nu^{-1}\|\q u^\ell\|_{\dot{B}_{2,1}^{-1+\frac {d}{2}} }\|\nu\q u^\ell\|_{\dot{B}_{2,1}^{1+\frac {d}{2}} }+\|V\|_{\dot{B}_{p,1}^{-1+\frac {d}{p}} }\|V\|_{\dot{B}_{p,1}^{1+\frac {d}{p}} }\Big)dt\nonumber\\
&\lesssim M+M^2+\delta^2+\nu^{-1}{\eta}+\nu^{-1}{\eta}^2.
\end{align}
Inserting \eqref{TT11} into \eqref{T8} and using \eqref{T1} imply that
\begin{align}\label{T9}
\z+\w
\le&\exp\Big(C(M+M^2+\delta^2+\nu^{-1}{\eta}+\nu^{-1}{\eta}^2)\Big)\Big\{\nu^{-\frac 12}{\eta}M\nonumber\\
&+\nu^{-1} {\eta}(M+{\eta}+\delta)\w+\nu^{-1}(M+{\eta}+\delta){\eta}+\nu^{-1} {\eta}M^2\Big\}.
\end{align}
Hence, assuming in addition that
 \begin{equation}\label{T10}
\nu^{-1}{\eta}+\nu^{-1}{\eta}^2\leq M+M^2\quad\hbox{and}\quad \delta\leq \max\{M,1\}.
 \end{equation}
We get from \eqref{T6} and \eqref{T9}  that
\begin{align}\label{T11}
\z+\w
\le&\exp(C(M+M^2))\Big\{\nu^{-\frac 12}{\eta}M+\nu^{-1} {\eta}(M+{\eta})\w+\nu^{-1}{\eta}^2+\nu^{-1} {\eta}M^2\Big\}
\end{align}
and
\begin{align}\label{T15}
\x+\y\le&\exp(C(M+M^2))
\Bigg\{\x(0)+M^2+1+\nu^{-1}(1+{\eta})({\eta}+M)\x+\nu^{-1}M\y\Bigg\}.
\end{align}
Therefore,  assume that
 \begin{equation}\label{T12}
 \nu^{-1}(1+{\eta})({\eta}+M)e^{CM} \ll 1.
  \end{equation}
Then we  get from \eqref{T11} and  \eqref{T15} that
\begin{align}
\z+\w
\le&{\eta}\exp(C(M+M^2))\Big\{\nu^{-\frac 12}M+\nu^{-1}({\eta}+ M^2+1)\Big\},\\
\x+\y\le&\exp(C(M+M^2))
(\x(0)+M^2+1).
\end{align}

Actually, it is natural to set
\begin{equation}\label{T18}
{\eta}= C\exp(C(M+M^2))(\x(0)+M^2+1)\quad\hbox{and}
\end{equation}
\begin{equation}\label{T19}
\delta=  C\exp(C(M+M^2))\bigl(\x(0)+M^2\bigr) \Bigl( \nu^{-1/2}M +\nu^{-1}(\x(0)+ M^2+1\bigr)\Bigr)\cdotp \end{equation}
Then take a suitably large constant $C$ such that
\begin{equation}\label{eq:nu}
C\exp(C(M+M^2))(\x(0)+1+M^2)\leq\sqrt{\nu},
\end{equation}
which implies the validity of the assumptions   \eqref{tiaojian1}, \eqref{T10} and  \eqref{T12}.

Consequently, due to   \eqref{eq:nu}, we define  $\eta$ and $\delta$ according to \eqref{T18} and \eqref{T19}
 so that \eqref{T1} is valid for all $T<T*.$ With the use of the global estimates of  $a$ and $v$, we conclude that $T*=+\infty$
and hence  \eqref{T1} is satisfied for all time.

The proof of Theorem \ref{dingli} is complete.

\bigskip
 \noindent \textbf{Acknowledgement.} This work is supported by the National Natural Science Foundation of China  under the grants 11601533 and 11571240.


\end{document}